\documentclass[12pt,a4paper]{article}
\usepackage{amsmath, amssymb}
\newtheorem{thm}{Theorem}
\newtheorem{prop}{Proposition}

\newtheorem{lem}{Lemma}
\newtheorem{dfn}{Definition}
\newtheorem{rmk}{Remark}

\begin{document}
\title{Projective normality of \\
 nonsingular toric varieties of dimension three \thanks{2002 Mathematics Subject Classification.
 Primary 14M25; Secondary 52B20}}
\author{Shoetsu OGATA\thanks{e-mail:  ogata{\char'100}math.tohoku.ac.jp}\\
Mathematical Institute Tohoku University\\ Sendai 980-8578, Japan}
\date{February, 2010}

\maketitle

%%%%%%%%%%%%%%%%%%% ‹r' %%%%%%%%%%%%%%%%%%%%%%%%%%%%%%
%\footnote{ %2000MSC
%2000 \textit{Mathematics Subject Classification}.
%Primary 14M25; Secondary 52B20.
%}
%\footnote{ %keywords
%\textit{Key words and phrases}. 
%Toric varieties, projectively normal.
%}
%\footnote{ %Thanks
%$^{*}$Partly supported by the Grant-in-Aid for Scientific Research (C),
%Japan Society for the Promotion of Science. 
%‰ÈŒ¤"ï•â•@Šî"ÕŒ¤‹†(A)'̏ꍇ
%}

\begin{abstract}
We show that if an ample line bundle  on a nonsingular toric 3-fold 
satisfies that its double adjoint bundle has no global sections, then it
 is  normally generated.  
 As an application, we give a proof of the normal generation of the anti-canonical bundle on a nonsingular toric Fano 4-fold. 
\end{abstract}

%%%%%%%%%%%%%
\section*{Introduction} %Introduction'ÉŒ©o'µ"ԍ†'ð•t'¯'éê‡'Í*'ð'Æ'é

  It is known that any ample line bundle  on a projective nonsingular toric
variety  is very ample (cf.  \cite[Corollary 2.15]{Od}).  
We call an  invertible sheaf on an algebraic variety a {\it line bundle}.
A line bundle $L$ on a projective variety is called {\it  normally generated}
 if the multiplication map $\Gamma(L)^{\otimes i} \to
\Gamma(L^{\otimes i})$ is surjective for all $i\ge1$.  
If an ample line bundle $L$ is normally generated, then we can easily
see that $L$ is very ample.   Furthermore, if the variety $X$ is normal,
then a normally generated ample line bundle $L$ defines the embedding
$\Phi_L: X \to \mathbb{P}(\Gamma(L))$ of $X$ as a {\it projectively normal}
variety.  

When we would ask questions about defining ideals of projective varieties we usually
assume that the varieties are projectively normal.  For example, 
Sturmfels \cite{S} asked whether  any  projective
nonsingular toric varieties embedded by normally generated ample
line bundles are defined by only quadrics (see also Cox \cite{C}).
Before giving any answer to such questions  we have to check 
 whether the variety be projectively normal, or the very ample
 line bundle on the variety be normally generated.  
 
We have few criteria of normal generation even on toric varieties.

Koelman \cite{K} showed that any ample line bundle on a toric
surface is normally generated.
Ewald and Wessels
\cite{EW} showed that for an ample line bundle $L$ on a projective toric
variety of dimension $n$, the twisted bundle $L^{\otimes i}$ is very ample for
$i\ge n-1$, and Nakagawa \cite{N} proved that $L^{\otimes i}$ is 
normally generated for $i\ge n-1$ (see also  \cite[Theorem 1]{NO}), more precisely
 that the multiplication map
 $$
 \Gamma(L)\otimes \Gamma(L^{\otimes i}) \longrightarrow \Gamma(L^{\otimes (i+1)})
 $$
 is surjective for $i\ge n-1$.
Ogata \cite{Og2} showed that a very ample line bundle on a certain class of 
projective toric 3-folds is normally generated.
This class consists of toric varieties which  are quotients of the projective
3-space $\mathbb{P}^3$ by  action of  finite abelian groups, and it contains
 weighted projective 3-spaces. 

In this paper we shall prove the following theorem 
(Theorem~\ref{sect4:tm} in Section~\ref{sect6}). 
   
 \begin{thm}\label{int:tm}
 Let $X$ be a nonsingular projective toric variety of dimension three.  Then
 any ample line bundle $L$ on $X$ satisfying that
 $H^0(X, L+2K_X)=0$ is normally generated.
 \end{thm}
 
 For a proof of Theorem~\ref{int:tm} we use the following result.
 \begin{thm}[Fakhruddin\cite{Fa}]\label{int:t2}
Let $X$ be a nonsingular projective toric surface.  Then, for an ample line bundle $A$ and a nef
line bundle $B$ on $X$, the multiplication map
$$
\Gamma(A)\otimes \Gamma(B) \longrightarrow \Gamma(A\otimes B)
$$
is surjective.
\end{thm}
 Kondo and Ogata \cite{KO}, and Haase, Nill, Paffenholz and Santos \cite{HN} generalized this
 to the case of singular toric surfaces. On the other hand,
 Ikeda \cite{Ik} generalized this to a certain class of nonsingular toric varieties in higher dimension.
 The class contains nonsingular
 toric varieties with the structure of $\mathbb{P}^r$-bundles over toric surfaces or over the projective
 line.  In particular, from \cite{Ik} we see that ample line bundles on a nonsingular toric variety
 which is a $\mathbb{P}^r$-bundle over a toric surface or over the projective line are
 normally generated.
 
  As an application of  Theorem~\ref{int:tm} we obtain a proof of the following.
  We do not use the classification of Fano polytopes (See Batyrev\cite{B}, or \O bro\cite{Ob}).
 \begin{thm}
Let $X$ be a nonsingular toric Fano variety of dimension four.  Then the anti-canonical
bundle $\mathcal{O}_X(-K_X)$ is normally generated.
\end{thm}

The structure of this paper is as follows:

In Sections 1 and 2 we recall basic results about toric varieties and line bundles on them.
A pair $(X, L)$ of a complete toric variety $X$ of dimension $n$ and an ample line bundle $L$ on $X$
corresponds to an integral convex polytope $P$ of dimension $n$ so that a basis of the space
of global sections of $L$ is parametrized by the set of all lattice points in $P$.
We recall that $L$ is normally generated if and only if the equalities
$$
(lP)\cap \mathbb{Z}^n+ P\cap \mathbb{Z}^n =((l+1)P) \cap \mathbb{Z}^n
$$
hold for all integers $l$.
In practice, it is enough to show the equalities for $1\le l\le n-2$.
In our case $n=3$, hence, $l=1$.  For a proof of Theorem~\ref{int:tm}, 
we need to investigate  properties of convex polytopes of dimension three.

In Section 3 we obtain a coarse classification of nonsingular integral convex polytopes
of dimension three without interior lattice points (Proposition~\ref{sect2:p1}).
By using the classification we prove a special case of Theorem~\ref{int:tm}, that is,
if the adjoint bundle of $L$ has no global sections, then $L$ is normally generated.
The condition on the adjoint bundle is equivalent that $P$ contains no lattice points
in its interior.  This is given by Proposition~\ref{sect2:p2}.

In Section 4 we treat the case that the adjoint bundle of $L$ has global sections
and investigate singularities of the interior polytope that corresponds to
the adjoint bundle of $L$ (Proposition~\ref{sect3:p2}).
For investigating the interior polytope,
 we can reduce to the case that the adjoint bundle is globally generated
(Proposition~\ref{sect3:p1}).

In Section 5 we treat the case that the interior polytope has adjacent singular vertices.
We show that if the interior polytope $Q$ has a singular vertex and if it has no lattice
points in its interior, then the adjoint bundle is normally generated (Propositions~\ref{sect5:p3}
and \ref{sect5:p4}).

In Section 6 we complete a proof of Theorem~\ref{int:tm}.

In Section 7 we give as an application a proof of the normal generation of the anti-canonical
bundle of a nonsingular toric Fano 4-fold (Proposition~\ref{sect7:p1}).
In our proof we do not use the classification of Fano polytopes.

%%%%%%%%%%%%%

\section{Projective toric varieties}\label{sect1}

In this section we recall the fact about toric varieties needed in this paper following Oda's
book \cite{Od}, or Fulton's book \cite{Fu}.
For simplicity, we consider toric varities are defined over the complex number field.

Let $N$ be a free $\mathbb{Z}$-module of rank $n$, $M$ its dual and 
$\langle, \rangle: M\times N \to\mathbb{Z}$ the canonical pairing.  
By scalar extension to the field $\mathbb{R}$ of real
numbers, we have real vector spaces $N_{\mathbb{R}}:=N\otimes_{\mathbb{Z}}\mathbb{R}$
and $M_{\mathbb{R}}:= M\otimes_{\mathbb{Z}}\mathbb{R}$.  We denote the same $\langle, \rangle$
as the pairing of $M_{\mathbb{R}}$
and $N_{\mathbb{R}}$ defined by scalar extension.
Let $T_N:= N \otimes_{\mathbb{Z}}\mathbb{C}^* \cong (\mathbb{C}^*)^n$ be the algebraic 
torus over the field $\mathbb{C}$ of complex numbers, where $\mathbb{C}^*$ is the multiplicative
group of $\mathbb{C}$.  Then $M=\mbox{Hom}_{\mbox{gr}}(T_N, \mathbb{C}^*)$ is the 
character group of $T_N$ and $T_N=\mbox{Spec}\ \mathbb{C}[M]$.
For $m\in M$ we denote $\bold{e}(m)$ as the character of $T_N$.  Let $\Delta$ be a finite complete
fan in $N$ consisting of strongly convex rational polyhedral cones $\sigma$ in $N_{\mathbb{R}}$,
that is, with a finite number of elements $v_1, \dots, v_s$ in $N$ we can write as
$$
\sigma = \mathbb{R}_{\ge0} v_1 + \cdots + \mathbb{R}_{\ge0} v_s
$$
and it satisfies that $\sigma \cap \{-\sigma\} =\{0\}$.
Then we have a complete toric variety $X=T_N \mbox{emb}(\Delta):= \cup_{\sigma \in \Delta}
U_{\sigma}$ of dimension $n$ (see  \cite[Section 1.2]{Od}, or  \cite[Section 1.4]{Fu}).
Here $U_{\sigma} =\mbox{Spec}\ \mathbb{C}[\sigma^{\vee}\cap M]$ and
$\sigma^{\vee}:=\{y\in M_{\mathbb{R}}; \langle y, x\rangle \ge0 \quad\mbox{for all $x\in \sigma$}\}$
is the dual cone of $\sigma$.
For the origin $\{0\} \in \Delta$, the affine open set $U_{\{0\}} = \mbox{Spec}\ \mathbb{C}[M]$
is the unique dense $T_N$-orbit.  We note that a toric variety 
defined by a fan is always normal.

If $|\Delta|:=\cup_{\sigma \in \Delta} \sigma =N_{\mathbb{R}}$, then the variety $X$ is complete.
Set $\Delta(s):=\{\sigma \in \Delta; \dim\sigma =s\}$.  Then $\tau \in \Delta(s)$ corresponds to
the $T_N$-orbit $\mbox{Spec}\ \mathbb{C}[\tau^{\perp}\cap M]$ and its closure $V(\tau)$, which
is also a $T_N$-invariant subvariety of dimension $n-s$.  Hence $\Delta(1)$ corresponds
to $T_N$-invariant irreducible divisors.  If  any cone $\sigma\in \Delta(n)$ of dimension $n$ 
is {\it nonsingular}, that is, there exist a 
$\mathbb{Z}$-basis $v_1, \dots, v_n$ in $N$ such that
\begin{equation}\label{sect1:eq1}
\sigma = \mathbb{R}_{\ge0} v_1 + \cdots + \mathbb{R}_{\ge0} v_n,
\end{equation}
then the toric variety $X$ is nonsingular.

Let $L$ be an ample $T_N$-equivariant line bundle on $X$.  Then
we have an integral convex polytope $P$ in $M_{\mathbb{R}}$ with
\begin{equation}\label{sect1:eq2}
H^0(X, L) \cong \bigoplus_{m\in P\cap M}\mathbb{C}\bold e(m),
\end{equation}
where $\bold e(m)$ are considered as rational functions on $X$ because
they are functions on an open dense subset $T_N$ of $X$ (see 
  \cite[Section 2.2]{Od}, or \cite[Section 3.5 ]{Fu}).  
Here an integral convex polytope $P$ in $M_{\mathbb{R}}$ is the convex
hull $\mbox{Conv}\{u_1, u_2, \dots, u_s\}$ in $M_{\mathbb{R}}$ of
a finite subset $\{u_1, u_2, \dots, u_s\} \subset M$.
We note that $\dim_{\mathbb{R}} P=\dim_{\mathbb{C}} X$.
The $l$ times twisted sheaf $L^{\otimes l}$ corresponds to 
the convex polytope
$lP :=\{lx\in M_{\mathbb{R}}; x\in P\}$.

On the other hand, for an integral convex polytope $P$ in $M_{\mathbb{R}}$ of dimension $n$
we can construct a projective toric variety $X$ of dimension $n$ and an ample line bundle
 $L$ satisfying (\ref{sect1:eq2}) (see  in \cite[Theorem 2.22]{Od}).
Indeed, for  each vertex $u_i$ of $P$ $(i=1, 2, \dots, r)$ we make
convex cone  $\mathbb{R}_{\ge0}(P - u_i):=\{\lambda(x-u_i)\in \mathbb{R}^n;
x\in P\quad \mbox{and} \quad \lambda \ge0\}$
and its dual cone $\tau_i$ in $N_{\mathbb{R}}$.   Set $\Delta$ to be a finite
complete fan of $N$ consisting of all faces of cones $\tau_i$ for
$i=1, 2, \dots, r$.  Then we obtain a projective toric variety
$X= T_N\mbox{emb}(\Delta)$ and an ample line bundle $L$
satisfying (\ref{sect1:eq2}).  In this sense we say that $P$ corresponds to
the {\it polarized toric variety} $(X, L)$.
If $X$ is nonsingular, then each cone $\tau_i$ has the same form as (\ref{sect1:eq1})
with a $\mathbb{Z}$-basis $v_1, \dots, v_n$ of $N$.  Hence the dual cone
$\tau_i^{\vee}= \mathbb{R}_{\ge0}(P-u_i)$ is a simplicial cone generated by
a $\mathbb{Z}$-basis $m_1-u_i, \dots, m_n-u_i$ of $M$.

\begin{dfn}
An integral convex polytope $P$ in $M_{\mathbb{R}}$ of dimension $n$ is
called {\it nonsingular} if for each vertex $u$ of $P$ the cone $\mathbb{R}_{\ge0}(P-u)$
is nonsingular in the sense of (\ref{sect1:eq1}).
We note that a nonsingular polytope $P$ is {\it simple}, that is, each vertex of $P$ is
contained in just $n$ faces of dimension $n-1$, or equivalently each vertex is contained
in just $n$ faces of dimension one.
\end{dfn}

We recall the notion that $L$ is {\it very ample}, that is, there is an
embedding of $X$ defined by the global sections of $L$:
$$
\Phi: X \to \mathbb{P}(H^0(X, L)).
$$
We can also interpret the condition for $L$ to be very ample  in terms of
$P$ as the condition that for each vertex $u$ of $P$ the semigroup
$\mathbb{R}_{\ge0}(P - u)\cap M$ in the cone $\mathbb{R}_{\ge0}(P - u)$
is generated by $(P - u)\cap M$.  
In other words, for each natural number $l$ all lattice points $x$ in
$l(P-u)$ are represented as a finite sum of elements $y_1, \dots, y_s$ in 
$(P-u)\cap M$.
We note that the number $s$ of elements $\{y_1, \dots, y_s\}$ in
$(P-u)\cap M$ needed for writing $x$ as their sum may be
different from $l$ such that $x$ lies in $l(P-u)$.
It is easy to see that any ample line bundle on a nonsingular toric variety
is very ample.

\begin{dfn}  An ample line bundle $L$ on a projective
variety $X$ is called {\it  normally generated} if the multiplication
map $\mbox{\textup{Sym}}^lH^0(X, L) \to H^0(X, L^{\otimes l})$ is
surjective for all $l\ge1$.  
\end{dfn}

\begin{dfn}
An integral convex polytope in $M_{\mathbb{R}}$ is called {\it normally generated} if
for the corresponding polarized toric variety $(X, L)$ the ample line bundle $L$ is
normally generated.
\end{dfn}

\begin{rmk}\label{rm}
If $X$ is toric and if $(X, L)$ corresponds with
an integral convex polytope $P$ in $M_{\mathbb{R}}$ satisfying (\ref{sect1:eq2}),
then the  normal generation of $L$ is equivalent to the condition
that for all $l\ge1$ every element $v\in lP\cap M$ be written as a sum
 $v=u_1 + \cdots + u_l$ of $l$ lattice points $u_i \in P\cap M$,
in other words, the condition that
$$
(lP)\cap M + P\cap M = ((l+1)P)\cap M \quad \mbox{for all $l\ge1$}.
$$
\end{rmk}

\section{Line bundles on toric varieties.}\label{sect2}

Let $\Delta$ be a complete fan of $N$ and let $X=T_N\mbox{emb}(\Delta)$
the corresponding toric variety.  For a cone of dimension one $\rho \in \Delta(1)$
we denote the primitive element of $\rho\cap N$ by $n(\rho)$.
Recall a $\Delta$-linear support function $h: N_{\mathbb{R}} \to \mathbb{R}$,
which is a continuous function linear on each cone $\sigma\in \Delta$,
defines a $T_N$-invariant Cartier divisor
$$
D_h := -\sum_{\rho \in \Delta(1)} h(n(\rho)) V(\rho)
$$
and an equivariant line bundle $\mathcal{O}_X(D_h)$
(see \cite[Section 2.1]{Od}).
For this line bundle from  \cite[Lemma 2.3]{Od}
we have an expression of the space of global sections as
\begin{equation}\label{sect1:eq3}
\Gamma(X, \mathcal{O}_X(D_h)) \cong \bigoplus_{m\in \Box_h} \mathbb{C}{\bold e}(m),
\end{equation}
where 
$$
\Box_h:= \{ m\in M_{\mathbb{R}}; \langle m,v\rangle \ge h(v)\quad\mbox{for all $v \in N_{\mathbb{R}}$}\}
$$
is a compact convex polytope in $M_{\mathbb{R}}$ (may be an empty set).
By definition there exist $l_{\sigma}\in M$ such that $h(v)= 
\langle l_{\sigma}, v\rangle$ for all $v\in \sigma$.
And we see that $\mathcal{O}_X(D_h)$ coincides with $\mathcal{O}_X\cdot {\bold e}(l_{\sigma})$
if they are restricted to $U_{\sigma}$.

A line bundle $L$ on $X$ is called {\it generated by global sections}, or shortly
{\it globally generated} if the map
$\Gamma(X, L)\otimes_{\mathbb{C}} \mathcal{O}_X \to L$ is surjective.

\begin{lem}[\cite{Od}]\label{sect1:l1}
For a complete toric variety $X=T_N\mbox{\rm emb}(\Delta)$ and a $\Delta$-linear support function
$h$ the following conditions are equivalent.
\begin{itemize}
 \item[{\rm (1)}] $\mathcal{O}_X(D_h)$ is globally generated.
 \item[{\rm (2)}] the linear system $|D_h|$ has no base points.
 \item[{\rm (3)}] $\Box_h = \mbox{\rm Conv}\{ l_{\sigma}; \sigma\in \Delta(n)\}$.
\end{itemize}
\end{lem}

From this and from the construction of polarized toric varieties 
we have a result of Mavlyutov \cite{Ma}.
\begin{lem}[Mavlyutov \cite{Ma}]\label{sect1:l2}
For a globally generated line bundle $\mathcal{O}_X(D_h)$ there exist an equivariant
surjective morphism $\pi: X \to Y$ to a toric variety $Y$ and an ample line bundle $A$
on $Y$ such that $\mathcal{O}_X(D_h)\cong \pi^*A$.
\end{lem}

From this lemma we see that $\mathcal{O}_X(D_h)$ is globally generated if and only if
$D_h$ is {\it nef} (see also  \cite[Theorem 3.1]{Ms}).

\section{Convex polytopes without interior lattice points}\label{sect3}

In this section we prove  Theorem~\ref{int:tm} in the case that 
$\Gamma(L
\otimes\mathcal{O}_X(K_X))=0$.

Let $X$ be a nonsingular projective toric 3-fold and $L$ an ample line bundle on $X$.
Let $P$ be the integral convex polytope of dimension three corresponding to
the polarized toric variety $(X, L)$.  From  \cite[Theorem 3.6]{Od} we have
\begin{equation}\label{sect2:eq1}
\Gamma(X, L\otimes \mathcal{O}_X(K_X)) \cong 
\bigoplus_{m\in \mbox{{\scriptsize Int}}(P)\cap M}
\mathbb{C}{\bold e}(m).
\end{equation}
Hence we see that $\Gamma(L\otimes\mathcal{O}_X(K_X))=0$ is equivalent to 
$\mbox{Int}(P)\cap M=\emptyset$.
In this section we consider an integral convex polytope $P$ of dimension three
satisfying the condition that $\mbox{Int}(P)\cap M=\emptyset$.

First we explain typical examples of nonsingular integral convex polytope
$P$ with $\mbox{Int}(P)\cap M=\emptyset$.
Set $P_0:=\mbox{Conv}\{(0,0,0), (1,0,0), (0,1,0), (0,0,1)\}$.  Then $P_0$
 defines the polarized toric variety $(\mathbb{P}^3, \mathcal{O}(1))$.
 Thus we see that $lP$ does not contain lattice points in its interior for $l=1,2,3$.
 We note that $lP_0$ is normally generated for all $l\ge1$.
 
 Set  $P_1:=\mbox{Conv}\{(0,0,0), (2,0,0), (0,2,0), (1,0,1), (0,1,1), (0,0,1)\}$ and
 $P_2^{(1)}:= \mbox{Conv}$ $\{ (0,0,0), (3,0,0), (0,3,0), (1,0,2), (0,1,2), (0,0,2)\}$.
 If you cut $2P_0$ at the hight 1, then you obtain $P_1$.
 The polytope $P_2^{(1)}$ is  cutt at the hight 2 from $3P_0$.
 Thus we see that $\mbox{Int}(P_1)\cap M=\mbox{Int}(P_2^{(1)})\cap M= \emptyset$.
The convex polytopes $P_1$ and $P_2^{(1)}$ define the blowing up of $\mathbb{P}^3$ 
at a $T_N$-invariant point.  This is also a toric $\mathbb{P}^1$-bundle over $\mathbb{P}^2$,
that is, $X\cong \mathbb{P}(\mathcal{O}\oplus \mathcal{O}(1))$.
We also have convex polytopes defining the blowing up of $\mathbb{P}^3$ at several points.
We may write $P_2^{(1)} = (3P_0)\cap \{0\le z\le2\}$.  Then we set
 $P_2^{(2)}:= P_2^{(1)}\cap \{0\le x\le2\}$, $P_2^{(3)}:= P_2^{(2)}\cap \{0\le y\le2\}$
and $P_2^{(4)}:=P_2^{(3)}\cap \{1\le x+y+z\le3\}$.  See the Figure~\ref{fig0} (a).
 
 \begin{figure}[h]
 \begin{center}
 \begin{tabular}{lr}
 \setlength{\unitlength}{1mm}
  \begin{picture}(50,60)(15,5)
   \put(20,20){\vector(4,-1){30}}
   \put(20,20){\vector(0,1){30}}
   \put(20,20){\vector(2,1){30}}
   \put(2,15){\makebox(20,10){$(0,0,0)$}}
   \put(53,10){\makebox(10,10)[bl]{$x$}}
   \put(15,50){\makebox(10,10)[bl]{$z$}}
   \put(51,32){\makebox(10,10)[bl]{$y$}}
   \put(41.5,14){\line(2,1){11}}
   \put(42,14){\line(0,1){11}}
   \put(28,41){\line(5,-2){19}}
   \put(28,41){\line(2,-5){3}}
   \put(20,37){\line(4,-1){11}}
   \put(20,37){\line(2,1){8}}
   \put(47,33){\line(2,-5){5.5}}
   \put(42,25){\line(2,-1){10}}
   \put(42,25){\line(-5,4){12}}
   \put(35,5){\makebox(10,10)[r]{$(2,0,0)$}}
   \put(55,14){\makebox(10,10)[]{$(2,1,0)$}}
   \put(45,22){\makebox(10,10)[l]{$(2,0,1)$}}
   \put(43,37){\makebox(10,10)[br]{$(0,3,0)$}}
   \put(42,25){\circle*{1}}
   \put(47,33){\circle*{1}}
   \put(25,40){\makebox(10,10)[l]{$(0,1,2)$}}
   \put(22,22){\makebox(10,10)[tl]{$(1,0,2)$}}
   \put(2,30){\makebox(20,10){$(0,0,2)$}}
   \put(20,20){\circle*{1}}
   \put(42,14.5){\circle*{1}}
   \put(20,37){\circle*{1}}
   \put(52.5,19.5){\circle*{1}}
   \put(30.5,34.5){\circle*{1}}
   \put(28,41){\circle*{1}}
  \end{picture}&
  
  \setlength{\unitlength}{1mm}
  \begin{picture}(50,60)(-5,5)
   \put(20,20){\vector(4,-1){30}}
   \put(20,20){\vector(0,1){30}}
   \put(20,20){\vector(2,1){30}}
   \put(2,15){\makebox(20,10){$(0,0,0)$}}
   \put(53,10){\makebox(10,10)[bl]{$x$}}
   \put(15,50){\makebox(10,10)[bl]{$z$}}
   \put(51,32){\makebox(10,10)[bl]{$y$}}
   \put(40,15){\line(0,1){30}}
   \put(35,27){\line(0,1){23}}
   \put(35,28){\line(1,-3){4.5}}
   \put(20,35){\line(1,1){15}}
   \put(20,35){\line(2,1){20}}
   \put(35,50){\line(1,-1){5}}
   \put(35,6){\makebox(10,10){$(1,0,0)$}}
   \put(40,20){\makebox(10,10){$(0,1,0)$}}
   \put(43,40){\makebox(10,10){$(1,0,a)$}}
   \put(30,47){\makebox(10,10){$(0,1,b)$}}
   \put(5,30){\makebox(10,10){$(0,0,c)$}}
   \put(20,20){\circle*{1}}
   \put(40,15){\circle*{1}}
   \put(35,27){\circle*{1}}
   \put(35,50){\circle*{1}}
   \put(20,35){\circle*{1}}
   \put(40,45){\circle*{1}}
  \end{picture}\\
  \mbox{(a): $P_2^{(2)}$} & \mbox{(b): $P_{a,b,c}$}
  \end{tabular}
 \end{center}
\caption{typical $P$ with $(\mbox{Int}P)\cap M=\emptyset$}
 \label{fig0}
\end{figure}

For $a \ge b\ge c\ge1$, set 
 $$
P_{a,b,c}:= \mbox{Conv}\{(0,0,0),(1,0,0),(0,1,0), (1, 0,a),(0,1,b),(0,0,c)\}.
$$
This is a prism with the basic triangle as  its section and its three edges have length $a, b$
and $c$.  See the Figure~\ref{fig0} (b).
The convex polytope
 $P_{a,b,c}$ defines a toric $\mathbb{P}^2$-bundle over $\mathbb{P}^1$, that is,
$X\cong \mathbb{P}(\mathcal{O}(a)\oplus \mathcal{O}(b)\oplus\mathcal{O}(c))$.

For convenience of explanation, next, we fix a notation of lattice points in $P$ near
a face of dimension two.  
We call a face of dimension two a {\it facet} and a face of dimension one an {\it edge}, simply.
Let $F_0$ be a facet of $P$.  Since $P$ is nonsingular, $F_0$ is
also nonsingular.  Denote $\{u_0, u_1, \dots, u_r\}$ the set of vertices of $F_0$.
Assume that $u_i$ is adjacent to $u_{i+1}$ for $i=0, 1, \dots, r$ (set $u_{r+1}=u_0$).
Take $m_1\in M$ on the edge $\overline{u_0 u_1}$ of $F_0$ and $m_2\in M$ on
$\overline{u_0 u_r}$ so that $\{m_1-u_0, m_2-u_0\}$ be a $\mathbb{Z}$-basis of
$(\mathbb{R}F_0)\cap M$.  
Since $P$ is nonsingular, we can take the lattice point $m_3 \in M\cap P$ on
the other edge meeting with $F_0$
at $u_0$  so that $\{m_1-u_0, m_2-u_0, m_3-u_0\}$
be a $\mathbb{Z}$-basis of $M$.  By using this basis we may identify $M$ with
$\mathbb{Z}^3$.  Let $(x,y,z)$ be the coordinates of $M_{\mathbb{R}}\cong \mathbb{R}^3$.
For each $u_i$ we can take the other edge of $P$ meeting with $F_0$ at $u_i$ and $w_i\in P\cap M$
on the edge with the coordinate $z=1$.  See the Figure~\ref{fig1}.

\begin{figure}[h]
 \begin{center}
 \setlength{\unitlength}{1mm}
  \begin{picture}(100,60)(0,5)
   \put(20,20){\vector(4,-1){30}}
   \put(20,20){\vector(0,1){30}}
   \put(20,20){\vector(2,1){30}}
   \put(15,20){\makebox(10,10)[bl]{$u_0$}}
   \put(53,10){\makebox(10,10)[bl]{$x$}}
   \put(15,50){\makebox(10,10)[bl]{$z$}}
   \put(50,35){\makebox(10,10)[bl]{$y$}}
   \put(40,15){\line(1,0){30}}
   \put(40,15){\line(0,1){10}}
   \put(70,15){\line(1,1){10}}
   \put(70,15){\line(0,1){10}}
   \put(40,30){\line(6,1){20}}
   \put(40,30){\line(0,1){10}}
   \put(30,5){\makebox(10,10)[r]{$u_1$}}
   \put(65,5){\makebox(10,10)[]{$u_2$}}
   \put(35,27){\makebox(10,10)[l]{$u_r$}}
   \put(35,25){\makebox(10,10)[br]{$w_1$}}
   \put(60,17){\makebox(10,10)[t]{$w_2$}}
   \put(30,35){\makebox(10,10)[tr]{$w_r$}}
   \put(50,20){\makebox(15,15)[bl]{$F_0$}}
   \put(50,40){\makebox(15,15){$P$}}
   \put(65,20){\makebox(20,20){$\ddots$}}
   \put(30,17){\circle*{1}}
   \put(20,30){\circle*{1}}
   \put(30,25){\circle*{1}}
   \put(25,10){\makebox(10,10)[l]{$m_1$}}
   \put(25,20){\makebox(10,10)[tl]{$m_2$}}
   \put(2,25){\makebox(20,10){$w_0=m_3$}}
   \put(20,20){\circle*{1}}
   \put(40,15){\circle*{1}}
   \put(70,15){\circle*{1}}
   \put(40,25){\circle*{1}}
   \put(70,25){\circle*{1}}
   \put(40,30){\circle*{1}}
   \put(40,40){\circle*{1}}
  \end{picture}
 \end{center}
\caption{$P$ and $F_0$ centered at $u_0$}
 \label{fig1}
\end{figure}

Now set $P(F_0):= \{0\le z\le1\}\cap P$ and $G:=\{z=1\}\cap P \subset P(F_0)\cap P$.  
Then $P(F_0)$ is an integral convex polytope with faces $F_0$ and $G$.
If $\dim G\le2$, then $G$ is a face of $P$.  When $\dim G=0$, that is, when $w_0=w_1=\dots =w_r$,
 we see that $r=2$ and $P=P_0$ since $P$ is nonsingular.
When $\dim G=1$, we see $r=3$ since $P$ is simple. 
In this case, we may assume $w_0=w_1$, then we see that $u_1=m_1$ and $u_2$ has the coordinate
of the form $(1,a,0)$ since $F_0$ is nonsingular.
If we write as $u_3=(0,b,0), w_2=w_3=(0,c,1)$, then we see that $P\cong P_{a,b,c}$.

We assume that $\dim G=2$.
If $G$ is a facet of $P$, then
all $w_i$'s are distinct since $P$ is simple.  
On the other hand, we note that
if all $w_i$'s are distinct, then $G$ has the same number of vertices
as that of $F_0$ and $G$ is nonsingular.  
Furthermore, $P(F_0)$ defines a toric 3-fold which is a toric $\mathbb{P}^1$-bundle over a toric surface
$Y$ defined by $F_0$.  In this case we can prove $P(F_0)$ is normally generated from
\cite[Theorem 2.5]{Ik}.

\begin{figure}[h]
 \begin{center}
 \setlength{\unitlength}{1mm}
  \begin{picture}(100,60)(0,5)
   \put(20,20){\vector(4,-1){30}}
   \put(20,20){\vector(0,1){30}}
   \put(20,20){\vector(2,1){30}}
   \put(15,18){\makebox(10,10)[bl]{$u_0$}}
   \put(53,10){\makebox(10,10)[bl]{$x$}}
   \put(15,50){\makebox(10,10)[bl]{$z$}}
   \put(50,35){\makebox(10,10)[bl]{$y$}}
   \put(30,17){\line(1,0){30}}
   \put(30,17){\line(-2,3){10}}
   \put(60,17){\line(1,1){8}}
   \put(60,17){\line(1,3){6}}
   \put(68,25){\line(-1,5){2}}
   \put(40,30){\line(6,1){15}}
   \put(40,30){\line(0,1){10}}
   \put(52,9){\makebox(10,10)[r]{$u_2$}}
   \put(65,15){\makebox(10,10)[]{$u_3$}}
   \put(35,27){\makebox(10,10)[l]{$u_r$}}
   \put(35,25){\makebox(10,10)[br]{$w_1$}}
  \put(60,30){\makebox(10,10)[t]{$w_2=w_3$}}
   \put(30,35){\makebox(10,10)[tr]{$w_r$}}
   \put(50,20){\makebox(15,15)[bl]{$F_0$}}
   \put(50,40){\makebox(15,15){$P$}}
   \put(30,17){\circle*{1}}
   \put(20,31){\circle*{1}}
   \put(30,25){\circle*{1}}
   \put(25,7){\makebox(10,10)[l]{$m_1=u_1$}}
   \put(26,20){\makebox(10,10)[tl]{$m_2$}}
   \put(2,25){\makebox(20,10){$w_0=w_1$}}
   \put(20,20){\circle*{1}}
   \put(66,35){\circle*{1}}
   \put(60,17){\circle*{1}}
   \put(68,25){\circle*{1}}
   \put(40,30){\circle*{1}}
   \put(40,40){\circle*{1}}
  \end{picture}
 \end{center}
\caption{$P$ and $F_0$ centered at $u_0$}
 \label{fig1e}
\end{figure}
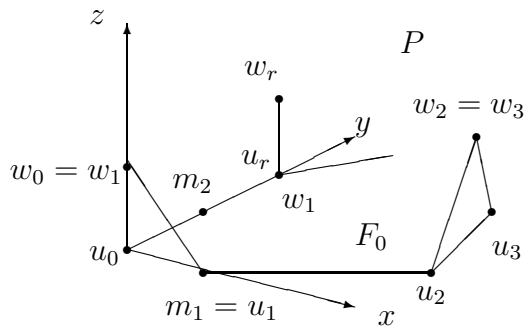

When $G=\{z=1\}\cap P$ is not a face of $P$, it may happen that $w_0=w_1$.
In this case, we see that $u_1=m_1$ because the facet $\mbox{Conv} \{u_0, u_1, w_0\}$
of $P$ is nonsingular.
If $w_0=w_1=w_2$, then $r=2$ and $u_2=m_2$, that is, $P=P_0$.
If $\dim G=2$ and if $w_0=w_1$, then $w_2\not=w_1$.  See the Figure~\ref{fig1e}.

In any case, since each edge $\overline{w_iw_{i+1}}$ of $G$ (with $w_i\not=w_{i+1}$) is
parallel to the edge $\overline{u_iu_{i+1}}$ of $F_0$, the integral polygon $G$ defines a nef
line bundle on $Y$.

\begin{prop}\label{sect2:p1}
Let $P$ be a nonsingular integral convex polytope in $M_{\mathbb{R}}$ of dimension three.
We assume that $P$ has no lattice points in its interior.  Then $P$ is one of the following.
\begin{itemize}
 \item[{\rm (1)}] $P$ is a convex hull of parallel two nonsingular facets 
 $F_0$ and $F_1$  such that
 the numbers of their vertices coincide.  This $P$ defines a toric $\mathbb{P}^1$-bundle
 over a nonsingular toric surface.
 \item[{\rm (2)}] $P$ is isomorphic to $P_0, 2P_0$, or $3P_0$.   The convex polytope $lP_0$
 corresponds to $(\mathbb{P}^3, \mathcal{O}(l))$.
 \item[{\rm (3)}] $P$ is isomorphic to $P_{a,b,c}$, or $2P_{a,b,c}$.  The convex polytope
 $P_{a,b,c}$ defines a toric $\mathbb{P}^2$-bundle over $\mathbb{P}^1$, that is,
 $\mathbb{P}(\mathcal{O}(a)\oplus \mathcal{O}(b)\oplus \mathcal{O}(c))$.
 \item[{\rm (4)}] $P$ is isomorphic to $P_2^{(i)}$ for $i=1, \dots, 4$.
 The convex polytope $P_2^{(i)}$ defines the blowing up of $\mathbb{P}^3$ at $T_N$-invariant
 $i$ points.
 In this case, we have $P=P(F_0)\cup P(F_1)$ by taking the parallel two facets $F_0$ and $F_1$.
 \item[{\rm (5)}] $P$ is contained in the triangular prism $\{(x,y,z)\in M_{\mathbb{R}};
 0\le x, 0\le y,  x+y\le2, 0\le z\}$. 
 In this case, we  have $P=P(F_1)\cup P(F_2)$ for a suitable choice of facets $F_i$.
\end{itemize} 
\end{prop}

{\it Proof.}  We use the notation described above.  
If $\dim G\le1$, then we see that $P\cong P_0$, or $P\cong P_{a,b,c}$ as discussed above.
In the following we assume $\dim G=2$.

Consider the case that $F_0$ and $G$ have the same number of edges.
If $G$ is a face of $P$, then it is in the case (1).

Assume that $G$ is not a facet of $P$.  Then the interior lattice
points $\mbox{Int}(G)\cap M$ are contained in the interior of $P$.
Thus by our assumption $G$ does not contain lattice points in its interior.

Set $G_0:=\mbox{Conv}\{(0,0), (1,0), (0,1)\}$ and $G_{a,b}:=\mbox{Conv}\{(0,0),  (0,1),$ 
$(a, 1),$ $(b,0)\}$ for $a\ge b\ge1$.  Then nonsingular integral convex polygons without interior
lattice points are only $G_0, 2G_0$ or $G_{a,b}$ up to affine transformations of $\mathbb{Z}^2$.
The convex polygons
$G_0$ and $2G_0$ correspond to the projective plane $\mathbb{P}^2$ with $\mathcal{O}(1)$
and $\mathcal{O}(2)$, respectively.  $G_{a,b}$ corresponds to the Hirzebruch surface
$\mathbb{P}(\mathcal{O}(a)\oplus \mathcal{O}(b))$ of degree $a-b$ with a suitable ample line
bundle.

If $G\cong G_0$, then we claim that $P\cong 2P_0$, or $P\cong P_{a,b,c}$.

\noindent
To see this, note that $F_0\cong kG_0$ for a positive integer $k$ since $F_0$ has to be a
nonsingular triangle.
In this case, $P(F_0)$ is combinatorially a prism, which may be given as $0\le z$, $0\le x$,
$0\le y$ and $x + y +(k-1)z \le k$.  The last three inequalities yield facets of $P$.
We distinguish two cases:

\noindent
If $k=1$, then $P=P_{a,b,c}$.

\noindent
If $k\ge2$, then the affine hyperplanes defined by the last three inequalities intersect in the point
$(0,0,\frac k{k-1})$, whose $z$-coordinate is less than or equal to 2, with equality only for $k=2$.
Since $G$ is not a facet of $P$, there has to exist a vertex of $P$ whose $z$-coordinate is greater than
or equal to 2.  Hence this implies $k=2$ and $P\cong 2P_0$.

If $G\cong 2G_0$, then it may happen $P\cong 3P_0, P\cong 2P_{a,b,c}$, or $P\cong P_2^{(1)}$.
In this case, we note that $F_0\cong kG_0$ for a positive integer $k$.
By the same reason above we have $k\le4$.  If $k=4$, then the point $(0,0,2)$ is a singular
vertex of the cone over $F_0$.  Hence, $1\le k\le3$.
We distinguish three cases:

\noindent
If $k=3$, then $P$ in contained in $3P_0$.
Set $F_1:=P\cap \{y=0\}$. Then $F_1$ is contained in the triangle $\mbox{Conv}\{(0,0,0),
(3,0,0), (0,0,3)\}$.
 If the point $(0,0,1)$ is a vertex of $F_1$, then there has to exist an edge
connecting $(0,0,1)$ with $(1,0,2)$ or $(2,0,1)$.  If the edge connects with $(1,0,2)$, then
the point $(1,0,2)$ is a singular vertex of $F_1$.  The situation is the same in the facet $P\cap
\{x=0\}$.  Thus, if the point $(0,0,1)$ is a vertex of $P$, then 
there have to exist two edges  connecting $(0,0,1)$ with $(2,0,1)$ and $(0,2,1)$,
hence, $G$ is a facet of $P$.  This contradicts to the assumption.  None of points
$(0,0,1), (2,0,1), (0,2,1)$ is a vertex of $P$.
If the point $(0,0,2)$ is a vertex of $P$, then it has to be connected with $(1,0,2)$ by an edge, hence, 
we have $P\cong P_2^{(1)}$, otherwise $P\cong 3P_0$.

\noindent
If $k=2$, then $P$ is contained in a prism with the twice of the basic triangle as its section,
which  may be given as $0\le z$, $0\le x$, $0\le y$ and $x+y\le2$.
It may happen $P\cong 2P_{a,b,c}$.  Even if not, this is the cese (5).

\noindent
If $k=1$, then we claim that $P$ is of the form (5), or $P\cong P_2^{(i)}$ for $i=1, \dots, 4$ (the case (4)).

\noindent
We assume that $P$ is not of the form (5).
Set $G':=P\cap\{z=2\}$.  Then $G'$ is a rational polygon.
We will prove that $G'$ contains the point $(1,1,2)$ as its interior.
We note that $G'$ is contained in the triangle 
$\tilde G:=\{0\le x, 0\le y, x+y\le3, z=2\}\cong 3G_0$.
The point $(1,1,2)$ is the center of $\tilde G$.  $G'$ is obtained by several cuts from $\tilde G$.

\begin{figure}[h]
\begin{center}
\begin{tabular}{lr}
  \setlength{\unitlength}{1mm}
  \begin{picture}(50,40)(5,5)
   \put(10,10){\line(1,0){30}}
   \put(10,10){\line(0,1){30}}
   \put(40,10){\line(-1,1){30}}
   \put(15,10){\line(1,-1){5}}
   \put(15,10){\line(-1,1){8}}
   \put(10,10){\circle*{1}}
   \put(10,20){\circle*{1}}
   \put(20,10){\circle*{1}}
   \put(20,20){\circle*{1}}
    \put(20,20){\circle*{1}}
    \put(20,30){\circle*{1}}
     \put(10,30){\circle*{1}}
     \put(10,40){\circle*{1}}
     \put(30,20){\circle*{1}}
     \put(40,10){\circle*{1}}
     \put(30,10){\circle*{1}}
   \put(-1,2){\makebox(10,10){$(0,0,2)$}}
   \put(13,7){\makebox(20,10){$(1,0,2)$}}
   \put(-2,17){\makebox(10,10){$(0,1,2)$}}
   \put(-2,33){\makebox(10,10){$(0,3,2)$}}
   \put(31,2){\makebox(30,10){$(3,0,2)$}}
   \put(14,18){\makebox(10,10){$(1,1,2)$}}
   \put(10,24){\makebox(10,10,){$G'$}}
    \put(25,24){\makebox(10,10){$\tilde G$}}
    \end{picture}&
    \setlength{\unitlength}{1mm}
  \begin{picture}(50,40)(-5,5)
   \put(10,10){\line(1,0){30}}
   \put(10,10){\line(0,1){30}}
   \put(40,10){\line(-1,1){30}}
   \put(20,15){\line(1,0){22}}
   \put(20,15){\line(-1,0){15}}
   \put(10,10){\circle*{1}}
   \put(10,20){\circle*{1}}
   \put(20,10){\circle*{1}}
   \put(20,20){\circle*{1}}
    \put(20,20){\circle*{1}}
    \put(20,30){\circle*{1}}
     \put(10,30){\circle*{1}}
     \put(10,40){\circle*{1}}
     \put(30,20){\circle*{1}}
     \put(40,10){\circle*{1}}
     \put(30,10){\circle*{1}}
   \put(-1,2){\makebox(10,10){$(0,0,2)$}}
   \put(14,2){\makebox(10,10){$(1,0,2)$}}
   \put(-2,17){\makebox(10,10){$(0,1,2)$}}
   \put(-2,33){\makebox(10,10){$(0,3,2)$}}
   \put(31,2){\makebox(30,10){$(3,0,2)$}}
   \put(14,18){\makebox(10,10){$(1,1,2)$}}
   \put(10,24){\makebox(10,10,){$G'$}}
    \put(25,24){\makebox(10,10){$\tilde G$}}
    \end{picture}\\
    \mbox{(a): cut at $(a-1)x+(b-1)y=1$} & \mbox{(b): cut at $(b-1)y=1$}
    \end{tabular}
   
   \end{center}
   \caption{$\tilde G$ containing $G'$}
   \label{fig1f}
   \end{figure}
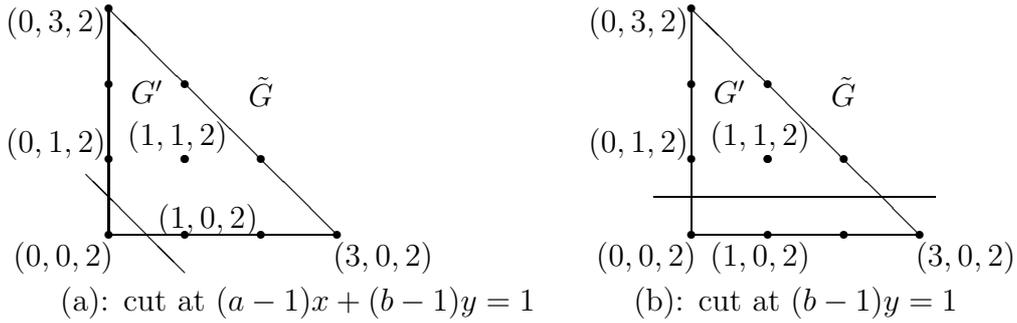

\noindent
If the point $(0,0,1)$ is a vertex of $P$, then there have to exist two edges connecting $(0,0,1)$
with $(1,0,a)$ and $(0,1,b)$.  
Since $G$ is not a facet of $P$, we see $a, b\ge1$ and one of them is greater than 1.
If both $a$ and $b$ are greater than 1, then $G'$ is obtained by  cut a triangle with the vertex
$(0,0,2)$ and with two edges of length $\frac1{a-1}$ $(\le1)$ and $\frac1{b-1}$ $(\le1)$.
See the Figure~\ref{fig1f} (a).
In this case, $G'$ contains $(1,0,2)$ and $(0,1,2)$, and the point $(1,1,2)$ is in the interior of $G'$.
If $a=1$ and $b=2$, then $G'$ is obtained by cut at the line $\{y=1\}$ and $G'\cong 2G_0$,
this implies that $P$ has the form (5) 
after an affine transform of $\mathbb{Z}^3$, since $G\cong 2G_0$.
  Hence, if $a=1$, then $b\ge3$.  In this case $G'$ is 
obtained by a cut at the line $\{y=\frac1{b-1}\}$
from $\tilde G$.  See the Figure~\ref{fig1f} (b).
The point $(1,1,2)$ remains in the interior of $G'$ since $b-1\ge2$.
Even if all three points $(0,0,1), (2,0,1)$ and $(0,2,1)$ are vertices of $P$, 
 then $(1,1,2)$ is the interior point of $G'$ unless $P$ is of the form (5).

\noindent
Since $P$ has no lattice points in its interior, $G'$ is a facet of $P$.  This implies 
that $a=b=2$.  This corresponds to $P_2^{(2)}$.  If $(2,0,1)$ or $(0,2,1)$ is
a vertex of $P$, then we have $P\cong P_2^{(3)}$ or $P\cong P_2^{(4)}$.

If $G\cong G_{a,b}$, then $F_0$ is a tetragon with two parallel edges.
 If $P$ is contained in the region $\{0\le y\le1\}$, then it is in the case (1).
If $F_0$ is a tetragon of the form $\mbox{Conv}\{(0,0), (0,k), (a',k), (b',0)\}$ with $k\ge2$,
then $k=2$ since $G$ is not a facet of $P$.
In this case, $P$ is contained in the prism $\{0\le x, 0\le y, 0\le z, x+z\le2\}$,
hence, we see that $P$ is of the form (5)
 by exchanging the role of $F_0$ with the facet of $P$
contained in the plane $\{x=0\}$.

Next we consider the case that $w_0=w_1$ in the Figure~\ref{fig1e}.
Then we see that $m_1=u_1=(1,0,0)$ and that $w_0$ is a vertex of $P$ since $P$ is nonsingular.
 If we write as
$u_2-u_1 =t(a,1,0)$, then  $a\ge -1$.   If $a=-1$, then $P\cong P_0$.  If $a\ge0$, then we can
reduce to the case treated above by exchanging the role of $F_0$ with the other face
$\mbox{Conv}\{u_0, u_1, w_0\}$.
\hfill $\Box$

We may apply the result of Ikeda \cite{Ik} to the cases (1) and (3) of Proposition~\ref{sect2:p1}
for the normal generation of $P$.  In this paper, we will use
 the following Lemmas for the normal generation of polytopes. 

\begin{lem}\label{sect2:l2}
Let $P$ be an integral convex polytope in $M_{\mathbb{R}}$.   If $P$ is a union of
normally generated integral convex polytopes, then $P$ is also normally generated.
\end{lem}
{\it Proof.}
Let $P=\cup_{i=1}^r Q_i$ be a decomposition into a union of integral convex polytopes
such that each $Q_i$ is normally generated.  For an integer $l$, take a lattice point in
$lP$, i.e., $m \in (lP)\cap M$.  Then we can choose $i$ so that $m \in lQ_i$ because
$lP=\cup_{i=1}^r lQ_i$.  Since $Q_i$ is normally generated, there exist $m_1, \dots,
m_l\in Q_i\cap M \subset P\cap M$ such that $m=m_1+\cdots +m_l$ from Remark~\ref{rm}.
\hfill $\Box$

\begin{lem}\label{sect2:l3}
The integral convex polytope $P(F_0)$ is normally generated.
\end{lem}
{\it Proof.}
We show that $P(F_0)\cap M +P(F_0)\cap M = (2P(F_0))\cap M$.
We note that $F_0$ and $G$ are normally generated because they are of dimension two.
From the result of Fakhruddin (Theorem~\ref{int:t2}), we see that
\begin{equation}\label{sect2:e3}
F_0\cap M + G\cap M = (F_0+G)\cap M
\end{equation}
because $F_0$ and $G$ define an ample and a nef line bundles on the nonsingular 
toric surface $Y$, respectively.

Take $m\in (2P(F_0))\cap M$.  If the $z$-coordinate of $m$ is 0, 1 or 2, then
$m$ is in $2F_0$, $F_0+G$ and $2G$, respectively.
Thus we can find $m_1, m_2 \in P(F_0)\cap M$ with $m=m_1+m_2$.
\hfill $\Box$

\begin{rmk}\label{sect2:rm1}
In the proof of lemma~\ref{sect2:l3} the equality~(\ref{sect2:e3}) is essential.
The result of Fakhruddin \cite{Fa} says that if each edge of $G$ has the same 
inner normal direction as that of some edge of $F_0$, then the equality~(\ref{sect2:e3})
holds.
The condition contains the case when 
$G$ is a line segment $E$ and $F_0$ is a tetragon with two edgs parallel to $E$.

Moreover, even if $F_0$ is a basic triangle, if the line segment $G=E$ is parallel to
an edge of the  basic triangle $F_0$, then the equality~(\ref{sect2:e3})
holds.
\end{rmk}

From Proposition~\ref{sect2:p1} and Lemmas~\ref{sect2:l2} and \ref{sect2:l3}
 we  prove Theorem~\ref{int:tm} of the special case.
\begin{prop}\label{sect2:p2}
Let $X$ be a projective nonsingular toric variety of dimension three and let $L$ 
an ample line bundle $L$ on $X$.
If $\Gamma(X, L\otimes\mathcal{O}_X(K_X))=0$, then $L$ is normally generated.
\end{prop}
{\it Proof.}
Let $P$ be the integral convex polytope corresponding to the polarized toric variety
$(X, L)$.  By the assumption $\Gamma(X, L\otimes\mathcal{O}_X(K_X))=0$,
the polytope $P$ does not contain lattice points in its interior.
We have a coarse classification  of such polytopes in Proposition~\ref{sect2:p1}.

We can apply Lemmas~\ref{sect2:l2} and \ref{sect2:l3} to 
the cases (1), (4) and (5) of Proposition~\ref{sect2:p1} for the normal generation of $P$.

If $P=kQ$ for some integral convex polytope $Q$ and $k\ge2$, then $P$ is normally generated from
\cite{N}.  The basic 3-simplex $P_0$ is trivially normally generated.
If $P\cong P_{a,b,c}$, then we can apply the result of Ikeda \cite{Ik}
since $P_{a,b,c}$ corresponds to a toric $\mathbb{P}^2$-bundle over the projective line,
or we may obtain the normal generation of $P_{a,b,c}$ from Remark~\ref{sect2:rm1}.
This completes the proof.
\hfill $\Box$

\section{Adjoint bundles.}\label{sect4}

In this section we investigate properties of the adjoint bundle 
$L\otimes \mathcal{O}_X(K_X)$
to an ample line bundle $L$ on $X$.

Let $L$ be an ample line bundle on a nonsingular projective toric variety
$X=T_N\mbox{emb}(\Delta)$ of dimension $n$.  Then there exists a $\Delta$-linear support
function $h: N_{\mathbb{R}} \to \mathbb{R}$ such that $L\cong \mathcal{O}_X(D_h)$ in the sense of Section~\ref{sect2}.  
Since ample line bundles on a toric variety are always globally generated,
that is, generated by its global sections, we have
$\Box_h =P$ from Lemma~\ref{sect1:l1} (3).  
Let $D=\sum_i D_i$ be the divisor consisting of all $T_N$-invariant irreducible
divisors on $X$.  Then $\mathcal{O}_X(K_X)\cong \mathcal{O}_X(-D)$.  We also have the
$\Delta$-linear support function $k: N_{\mathbb{R}} \to \mathbb{R}$ such that $\mathcal{O}_X
(K_X)\cong \mathcal{O}_X(D_k)$.  Here $k$ is defined on an $n$-dimensional cone
$\sigma =\sum_{i=1}^n \rho_i \in \Delta(n)$ by $k(n(\rho_i))=1$, where $\rho_i =\mathbb{R}_{\ge0}
n(\rho_i)$ and $n(\rho_i) \in \rho_i\cap M$ is the primitive element.
We want to describe the $\Delta$-linear support function $h+k$ of $L\otimes \mathcal{O}_X(K_X)$.

Assume that $\Gamma(L\otimes \mathcal{O}_X(K_X))\not=0$, equivalently that
$\mbox{Int}(P)\cap M \not=\emptyset$.  We know $\Box_{h+k} \cap M=\mbox{Int}(P)\cap M$.
Set $Q:=\mbox{Conv}(\mbox{Int}(P)\cap M)$.
We call $Q$ the {\it interior polytope} of $P$.  We see that $Q\subset \Box_{h+k}$
because $\Box_{h+k}$ is convex.

Let $u_0 \in P$ be a vertex of $P$.  Then there is the $n$-dimensional cone $\sigma \in \Delta(n)$
such that $\sigma^{\vee} \cong \mathbb{R}_{\ge0} (P-u_0)$.  We see that $u_0=l_{\sigma}$
in the sense of Section 2.  Since $\sigma$ is nonsingular, there are $m_1, \dots, m_n \in P\cap M$
such that $\{m_1-u_0, \dots, m_n-u_0\}$ is a $\mathbb{Z}$-basis of $M\cong \mathbb{Z}^n$ and that
$\mathbb{R}_{\ge0}(P-u_0) = \sum_{i=1}^n \mathbb{R}_{\ge0}(m_i-u_0)$.
Set $\bar{l}_{\sigma}:=u_0 +\sum_{i=1}^n (m_i-u_0)$.
Then the lattice point $\bar{l}_{\sigma}-u_0=\sum_{i=1}^n (m_i-u_0)$ 
is in the interior of $\sigma^{\vee}=\mathbb{R}_{\ge0}(P-u_0)$ 
and $(\mbox{Int}\ \sigma^{\vee})\cap M = (\bar{l}_{\sigma}-u_0) +\sigma^{\vee}\cap M$.
Thus we have $(h+k)(v)=\langle\bar{l}_{\sigma},v\rangle$ for all $v\in \sigma$.
Since the set of vertices of $P$ corresponds to $\Delta(n)$, we can define $\bar{l}_{\sigma}
\in M$ for all $\Delta(n)$.
If all $\bar{l}_{\sigma}$ are contained in $\Box_{h+k}$, then $L\otimes \mathcal{O}_X(K_X)$
is generated by global sections from Lemma~\ref{sect1:l1} and $\Box_{h+k}=Q$.
Unfortunately, $L\otimes \mathcal{O}_X(K_X)$ is not always generated by global sections.
Even if not, we will see $\Box_{h+k}=Q$ when $\dim X=3$ in the following Proposition.

The first statement of the Proposition is a corollary of the result of Fujita \cite[Theorems 1, 2 and 3]{Fj}.
We have to investigate the shape of $P$ when $D_h+K_X$ is not nef.

\begin{prop}\label{sect3:p1}
Let $X$ be a projective nonsingular toric variety of dimension three and let $L$ 
an ample line bundle $L$ on $X$.
If $\Gamma(X, L\otimes\mathcal{O}_X(K_X))\not=0$, then
there exists a polarized toric variety $(Y, A)$ of dimension three such that
$A\otimes \mathcal{O}_Y(K_Y)$ is globally generated and that 
the injective homomorphism $L \hookrightarrow \pi^*A$
 induces the isomorphism
$\Gamma(X, L\otimes\mathcal{O}_X(K_X))\cong \Gamma(Y, A\otimes \mathcal{O}_Y(K_Y))$,
 where $Y$ is a nonsingular toric variety obtained by contraction $\pi: X \to Y$ 
 of divisors to points.
 
 Moreover, we have $\Box_{h+k}=Q$.
 \end{prop}
{\it Proof.}
Let $u_0\in P$ be a vertex and  $F_0$ a facet  containing $u_0$.
The two edges of $F_0$ meeting at $u_0$ have the lattice points $m_1$ and $m_2$
respectively so that  $\{m_1-u_0, m_2 -u_0\}$ is a $\mathbb{Z}$-basis of
$(\mathbb{R}F_0) \cap M\cong \mathbb{Z}^2$.  Then we have the same
figure as the Figure~\ref{fig1} and the coordinate system $(x,y,z)$ of $M\cong \mathbb{Z}^3$. 

Consider the point $(1,1,1)$, which is $\bar{l}_{\sigma}$ of $\sigma^{\vee} =\mathbb{R}_{\ge0}
(P-u_0)$ as described above.
If $(1,1,1)$ is an interior lattice point of $P-u_0$, that is, if it is a vertex of $Q-u_0$,
 then $(1,1,1)$ is also the vertex of $\Box_{h+k}-u_0$. 
 
We assume that the point $(1,1,1)$ is not contained in $Q-u_0$.  
As in the proof of Proposition~\ref{sect2:p1}, we set $G:=(P-u_0)\cap \{z=1\}$.
We note that $G$ is not a facet of $P$ since $P$ contains interior lattice points.
Then the assumption implies that $(1,1,1)$ is not contained in the interior of $G$.
We may assume that $G$ contains the points $(1,0,1)$ and $(0,1,1)$.
If $(1,0,1)$ is not contained in $G$, then $(0,0,1)$ is a vertex of $P$ and the facet
$P\cap\{y=0\}$ is the basic triangle $\mbox{Conv}\{(0,0,0), (1,0,0), (0,0,1)\}$.  
In this case, if we  exchange the role of $F_0$
with the facet $P\cap\{y=0\}$, then new $G$ satisfies the assumption.

If $G$ is a tetragon, then it has two parallel edges with the distance one, hence, $F_0$ also
has two parallel edges.  In this case, $P$ has the shape of (1) or (5) in Proposition~\ref{sect2:p1}.
Then $P$ cannot contain lattice points in its interior.

If $G$ is a triangle not containing $(1,1,1)$, then $G$ has $(1,0,1)$ and $(0,a,1)$ ($a\ge1$) 
as its vertices.   
If $a=1$, then $P\cong P_{a,b,c}$ or $P$ is contained in $2P_0$ since $F_0$ has an edge 
parallel to the edge of $G$ connecting $(1,0,1)$ and $(0,1,1)$.
In both cases, $P$ cannot contain interior lattice points.
If $a\ge2$, then $(1,0,1)$ is a singular vertex of $G$, hence, it happens $w_1=w_2$ in
the Figure~\ref{fig1e}.  Since $F_0$ is nonsingular and has an edge parallel to 
the edge of $G$ connecting $(1,0,1)$ and $(0,1,a)$, we see that $w_2 - w_1$ has 
the direction $(0,1,0)$, hence, $F_0$ has two parallel edges. 
In this case, $P$ has the shape of (1) or (5) in Proposition~\ref{sect2:p1}.

\begin{figure}[h]
\begin{center}
 \begin{tabular}{lr}
  \setlength{\unitlength}{1mm}
  \begin{picture}(50,40)(0,5)
     \put(20,10){\line(1,0){20}}
   \put(20,10){\circle*{1}}
   \put(20,10){\line(-1,1){10}}
   \put(20,10){\line(-1,2){3.5}}
   \put(10,20){\line(0,1){20}}
   \put(10,20){\circle*{1}}
   \put(10,30){\circle*{1}}
   \put(10,20){\line(2,-1){6}}
   \put(17,17){\line(2,1){20}}
   \put(17,17){\circle*{1}}
   \put(10,10){\makebox(10,10)[bl]{$F_0$}}
   \put(20,5){\makebox(10,10)[bl]{$m_1$}}
   \put(4,19){\makebox(10,10)[bl]{$u_0$}}
   \put(3,30){\makebox(10,10)[bl]{$m_3$}}
   \put(18,10){\makebox(10,10){$m_2$}}
   \put(30,30){\makebox(10,10){$P$}}
   \put(25,23){\circle*{1}}
   \put(20,20){\makebox(10,10){$m_0$}}
  \end{picture} &

  \setlength{\unitlength}{1mm}
  \begin{picture}(50,40)(0,0)
   \put(10,10){\line(1,0){30}}
   \put(10,10){\line(1,1){20}}
   \put(10,10){\line(2,1){20}}
   \put(10,10){\line(0,1){20}}
   \put(10,10){\circle*{1}}
   \put(25,10){\circle*{1}}
   \put(20,15){\circle*{1}}
   \put(10,20){\circle*{1}}
   \put(5,5){\makebox(10,10)[bl]{$O$}}
   \put(25,5){\makebox(10,10)[bl]{$v_1$}}
   \put(25,10){\makebox(10,10)[l]{$v_2$}}
   \put(5,15){\makebox(10,10)[l]{$v_3$}}
   \put(30,30){\makebox(10,10)[b]{$\rho_0$}}
   \put(25,10){\line(-3,2){15}}
   \put(25,10){\line(-1,1){5}}
   \put(10,20){\line(2,-1){10}}
  \end{picture}\\
  \mbox{(a): $P$ around $F_0$} & \mbox{(b): $\Delta$ around $\rho_0$}
 \end{tabular}
 \end{center}
 \caption{Local shapes of $P$ and $\Delta$}
 \label{fig2}
\end{figure}
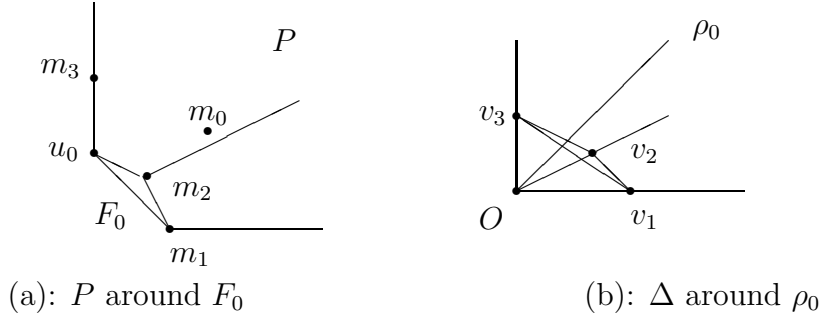

From this argument we see that $G$ is the triangle $\mbox{Conv}\{(0,0,1), (2,0,1), (0,2,1)\}\cong 2G_0$
containing the point $(1,1,1)$ in its boundary.  It is also contained 
in the boundary of $P-u_0$.  
As in the proof of Proposition~\ref{sect2:p1} (treating the case (4)), we see 
 that $F_0\cong G_0$ 
and that $(1,1,2)$ is an interior lattice point of $P-u_0$ since  $\mbox{Int}(P)\cap M\not=\emptyset$.
We note that $(1,1,2)$ is a vertex of $Q-u_0$.
We denote $(1,1,2)=m_0-u_0$ in $P -u_0$.  Then we see $m_0 \in Q \subset \Box_{h+k}$.
By taking an affine transformation of
$M\cong\mathbb{Z}^3$, we may set $u_0=(-1,-1,0), m_1=(0,-1,-1), m_2=(-1,0,-1),
m_3=(-1,-1,1)$.   Then the point $(1,1,2)$ in 
$P-u_0$ is transformed to the origin $m_0$.  See Figure~\ref{fig2}~(a).

The facet $F_0=\mbox{Conv}\{u_0, m_1, m_2\}$ corresponds to $(\mathbb{P}^2, \mathcal{O}
(1))$, which is a $T_N$-invariant divisor $V(\rho_0)$ on $X$ with $\rho_0 \in \Delta(1)$.
From the Figure~\ref{fig2}~(a) we may draw the picture of $\Delta$ around $\rho_0$.
See Figure~\ref{fig2}~(b).

Here $\{v_1, v_2, v_3\}$ is a $\mathbb{Z}$-basis of $N\cong \mathbb{Z}^3$ and the primitive
element of $\rho_0\cap N$ is $n(\rho_0) =v_1+v_2+v_3$.  In other words, $\rho_0$ gives
the barycentric subdivision of the nonsingular cone $\sum_{i=1}^3 \mathbb{R}_{\ge0}v_i$.
In terms of algebraic geometry, this is locally isomorphic to the blow up $\mathbb{C}^3$
at the origin.  We denote this blow up by $\pi: X \to Y$ and $E$ the exceptional divisor.
If we add the simplex $\mbox{Conv}\{(-1,-1,-1), m_1, m_2, u_0\}$ to $P$
in Figure~\ref{fig2}~(a), then we obtain an ample line bundle $A$ on $Y$
so that $L\otimes \mathcal{O}_X(E)\cong\pi^*A$.  Since $K_X=\pi^*K_Y +2E$, we have
$L\otimes\mathcal{O}_X(K_X) \cong \pi^*(A \otimes \mathcal{O}_Y(K_Y))\otimes
\mathcal{O}_X(E)$.   

Here we may write as
$$
\Box_h =\{m\in M_{\mathbb{R}}; \langle m, n(\rho)\rangle \ge h(n(\rho)) \quad
\mbox{for} \quad \rho \in \Delta(1)\}.
$$
Near $m_0$ we have
$$
m_0\in\Box_{h+k} \subset \cap_{i=1}^3 \{m\in M_{\mathbb{R}}; 
\langle m, v_i\rangle \ge h(v_i) +1 \},
$$
where $m_0$ is the apex of the triangular cone in the right hand side.
Since $\Box_{h+k}$ is an intersection of half-spaces, the point $m_0$ is a vertex of
$\Box_{h+k}$.  This implies $\Box_{h+k}=Q$.
\hfill $\Box$

Set $P'$ the integral polytope corresponding to $(Y, A)$.  Then Proposition~\ref{sect3:p1}
 implies that
both $P$ and $P'$ have the same interior polytope $Q$ and that $Q$ defines the nef
line bundle $A\otimes \mathcal{O}_Y(K_Y)$ on $Y$.

Set $\{F_i \subset P; i\in I\}$ and $\{F'_j\subset P';  j\in J\}$ the sets of all facets of $P$
and $P'$, respectively.  Then we have decompositions
$$
P =\cup_{i\in I} P(F_i) \cup Q \quad\mbox{and}\quad P' = \cup_{j\in J} P'(F'_j) \cup Q.
$$
In order to investigate the shape of $Q$, we may assume that $L\otimes\mathcal{O}_X(K_X)$
is globally generated.

In the following, we assume that $H^0(X, L\otimes\mathcal{O}_X(K_X))\not=0$,
 $L\otimes\mathcal{O}_X(K_X)$ is globally generated and that $\dim\ Q=3$.

Let $m_0$ be a vertex of $Q$ and $E_0 \subset Q$ a facet  containing $m_0$.
Then we can choose a facet $F_0 \subset P$  with the vertex $u_0$ such that
the primitive elements $m_1, m_2\in F_0\cap M$ and $m_3 \in P\cap M$ on three edges
meeting at $u_0$ form a $\mathbb{Z}$-basis of $M+u_0$ and that 
$m_0-u_0 =(m_1-u_0)+(m_2-u_0)+(m_3-u_0)$.  Moreover we may take $F_0$ so that
$G\cap Q= E_0$ as in the Figure~\ref{fig1} because $Q$ is surrounded by planes
parallel to $F_0$'s.
If $w_i\not=w_{i+1}$ for all $i$, then $G$ is a nonsingular integral polygon.  
It may happen that $w_0=w_1$ as in the Figure~\ref{fig1e}.

\begin{figure}[h]
\begin{center}
 \begin{tabular}{lr}
  \setlength{\unitlength}{1mm}
  \begin{picture}(50,30)(5,5)
   \put(10,10){\line(1,0){10}}
   \put(10,10){\line(0,1){20}}
   \put(20,10){\line(2,1){25}}
   \put(10,10){\circle*{1}}
   \put(10,20){\circle*{1}}
   \put(20,10){\circle*{1}}
   \put(40,20){\circle*{1}}
   \put(5,6){\makebox(10,10)[bl]{$u_0$}}
   \put(17,5){\makebox(20,10)[bl]{$u_1=m_1$}}
   \put(3,20){\makebox(10,10)[bl]{$m_2$}}
   \put(35,13){\makebox(30,10)[bl]{$(a+1,1,0)$}}
    \put(20,22){\makebox(10,10){$F_0$}}
    \end{picture} &

  \setlength{\unitlength}{1mm}
  \begin{picture}(50,30)(-5,5)
   \put(10,10){\line(0,1){20}}
   \put(10,10){\line(2,1){30}}
   \put(10,10){\circle*{1}}
   \put(20,20){\circle*{1}}
   \put(10,20){\circle*{1}}
   \put(30,20){\circle*{1}}
   \put(4,29){\makebox(10,10)[bl]{$w_r$}}
   \put(40,21){\makebox(10,10)[bl]{$w_2$}}
   \put(1,5){\makebox(20,10)[bl]{$w_0=w_1=(0,0,1)$}}
   \put(-4,20){\makebox(20,10)[bl]{$(0,1,1)$}}
   \put(30,13){\makebox(30,10)[bl]{$(a,1,1)$}}
    \put(20,22){\makebox(10,10){$G$}}
     \end{picture}\\
  \mbox{(a): $F_0$ near $u_0$ and $u_1$} & \mbox{(b): $G$ near $w_0$}
 \end{tabular}
 \end{center}
 \caption{Local shapes of $F_0$ and $G$}
 \label{fig2e}
\end{figure}
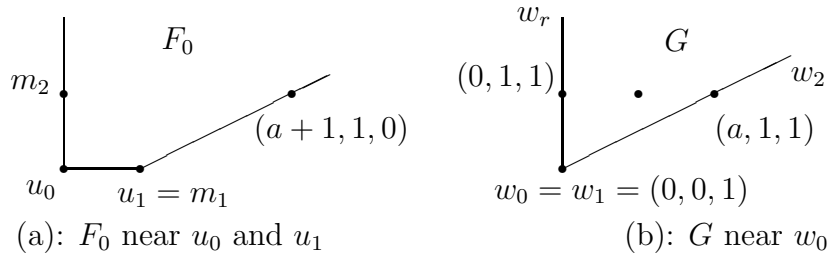

We have to investigate the shape of $G$ and the interior polygon of $G$.
If $w_0=w_1$ in the Figure~\ref{fig1e}, then we see $w_1\not= w_2$ and $w_0\not=w_r$.
In this case, the edge $\overline{u_1u_2}$ of $F_0$ has the direction $(a,1,0)$ with $a\ge1$
since $\{u_0-u_1, u_2-u_1\}$ is a part of $\mathbb{Z}$-basis of $M$.
Thus the edge $\overline{w_1w_2}$ of $G$ has the same direction $(a,1,0)$.
If $a=1$, then $w_0=w_1$ is a nonsingular vertex of $G$, otherwise it is singular.
See the Figure~\ref{fig2e}.

In terms of algebraic geometry, the algebraic surface corresponding to $G$ is locally
obtained by the contraction of the projective line with the self-intersection number
$-a$ to a point.  When $a=2$ we call the singularity $A_1$-singularity, which is also a
rational double point.

\begin{lem}\label{sect3:l1}
Let $G$ be an integral convex polygon of dimension two whose vertices have
singularities at worst described in the Figure~\ref{fig2e} (b).  
We assume that no two singular vertices are adjacent.
If the interior polygon $\mbox{\rm{Conv}}\{(\mbox{\rm{Int}}\ G)\cap \mathbb{Z}^2\}$
 is of dimension two, then
it has at worst $A_1$-singularities.

In particular, if $G$ is a nonsingular polygon, then the interior polygon is also nonsingular.
\end{lem}
{\it Proof.}
We denote $G^{\circ}:=\mbox{Conv}\{(\mbox{Int}\ G)\cap \mathbb{Z}^2\}$ the interior
polygon of $G$.  We assume $\dim G^{\circ}=2$.

First we prove that if $G$ is nonsingular, then  $G^{\circ}$ is also nonsingular.
In this case, we may consider as $G=F_0$  in the Figure~\ref{fig2e} (a).
Let $u_0$ be a vertex of $G$. The two edges meeting at $u_0$ have the lattice points
$m_1$ and $m_2$ respectively so that $\{m_1-u_0, m_2-u_0\}$ is a basis of
$\mathbb{Z}^2$.  Take $u_0$ the origin and the coordinates $(x, y)$ of $\mathbb{R}^2$
as in the Figure~\ref{fig2e} (a) 
such that $\{m_1-u_0, m_2-u_0\}$ is the basis.  Then the point $(1,1)$ is in $G^{\circ}$.

If both two points $(2,1)$ and $(1,2)$ are contained in $G^{\circ}$, then we see that 
the point $(1,1)$ is a nonsingular vertex of $G^{\circ}$ since $(2,1), (1,2) \in \partial G^{\circ}$.

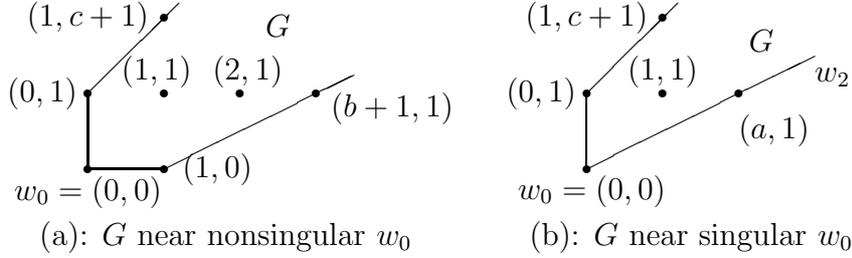
\begin{figure}[h]
\begin{center}
 \begin{tabular}{lr}
  \setlength{\unitlength}{1mm}
  \begin{picture}(50,30)(5,5)
   \put(10,10){\line(1,0){10}}
   \put(10,10){\line(0,1){10}}
   \put(20,10){\line(2,1){25}}
   \put(10,20){\line(1,1){12}}
   \put(10,10){\circle*{1}}
   \put(10,20){\circle*{1}}
   \put(20,10){\circle*{1}}
   \put(40,20){\circle*{1}}
   \put(20,20){\circle*{1}}
   \put(20,30){\circle*{1}}
   \put(30,20){\circle*{1}}
   \put(5,2){\makebox(10,10){$w_0=(0,0)$}}
   \put(17,5){\makebox(20,10){$(1,0)$}}
   \put(-1,15){\makebox(10,10){$(0,1)$}}
   \put(35,13){\makebox(30,10){$(b+1,1)$}}
   \put(14,18){\makebox(10,10){$(1,1)$}}
   \put(26,18){\makebox(10,10){$(2,1)$}}
   \put(5,25){\makebox(10,10){$(1,c+1)$}}
    \put(30,24){\makebox(10,10){$G$}}
    \end{picture} &

  \setlength{\unitlength}{1mm}
  \begin{picture}(50,30)(-5,5)
   \put(10,10){\line(0,1){10}}
   \put(10,10){\line(2,1){30}}
   \put(10,20){\line(1,1){12}}
   \put(10,10){\circle*{1}}
   \put(20,20){\circle*{1}}
   \put(10,20){\circle*{1}}
   \put(30,20){\circle*{1}}
   \put(20,30){\circle*{1}}
   \put(5,25){\makebox(10,10){$(1,c+1)$}}
   \put(40,21){\makebox(10,10)[bl]{$w_2$}}
   \put(1,5){\makebox(20,10)[bl]{$w_0=(0,0)$}}
   \put(-1,15){\makebox(10,10){$(0,1)$}}
   \put(30,13){\makebox(30,10)[bl]{$(a,1)$}}
   \put(15,18){\makebox(10,10){$(1,1)$}}
    \put(28,22){\makebox(10,10){$G$}}
     \end{picture}\\
  \mbox{(a): $G$ near nonsingular $w_0$} & \mbox{(b): $G$ near singular $w_0$}
 \end{tabular}
 \end{center}
 \caption{Local shapes of  $G$}
 \label{fig2f}
\end{figure}

%\noindent
If $(1,0)$ is not vertex or if it is a vertex of $G$ with the other edge going to
$(b+1,1)$ with $b\ge2$, then $G^{\circ}$ contains the point $(2,1)$ since $\dim G^{\circ}=2$.
See the Figure~\ref{fig2f} (a).
Since the situation at the point $(0,1)$ is the same, we set $(0,1)$ is a vertex with the other edge
of the direction $(1,1)$ and $(1,0)$ is a vertex with the other edge of the direction $(b,1)$
with $b\ge2$. 
In this case,  the points
$(2,1)$ and $(2,2)$ are contained in $G^{\circ}$ and
$(2,2)$ is in the boundary of $G^{\circ}$ since lattice points in $G$ are exhausted by
the lines $y=x+k$ with $k\le1$.
Thus, $(1,1)$ is a nonsingular vertex of $G^{\circ}$.

Next we set $w_0$ the singular vertex of $G$ as in the Figure~\ref{fig2f} (b).
In this case, if the point $(0,1)$ is a nonsingular 
vertex of $G$ and the edge from it connects with the point $(1, c+1)$,
then $(1,2) \in G^{\circ}$ for $c\ge2$.  When $c=1$ the point $(2,2)$ is contained in
 $G^{\circ}$ as a boundary point since lattice points in $G$ are exhausted by the lines
 $y=x+k$ with $k\le1$.
Thus we see that the point $(1,2)$ or $(2,2)$ is on the edge of $G^{\circ}$.

%\noindent
If $a\ge2$, then the points $(1,1)$ and $(2,1)$ are contained in $G^{\circ}$.
In this case, the point $(1,1)$ is a nonsingular vertex of $G^{\circ}$.

%\noindent
Assume $a=1$.  We distinguish two cases:

\noindent
If the point $(2,1)$ is not vertex or if it is a vertex and the edge from it connects with 
the point $(2\alpha +1,
\alpha+1)$ with $\alpha\ge2$, then the point $(3,2)$ is contained in $G^{\circ}$
as a boundary point since lattice points in $G$ are exhausted by the lines
$2y=x+k$ with $k\ge0$.
In this case, the vertex $(1,1)$ of $G^{\circ}$ is nonsingular if $(0,1)$ is a vertex with
the other edge going to $(1,2)$ or it is $A_1$-singularity otherwise.

\noindent
If $\alpha=1$, then $(3,2)\in \partial G$ and $G^{\circ}$ contains $(2,2)$
as its boundary point since lattice points in $G$ are exhausted by the lines
$y=x+k$ with $k\ge -1$.
In this case we see that the point $(1,1)$ is a nonsingular vertex of $G^{\circ}$
 since $\dim G^{\circ}=2$.

\begin{figure}[h]
\begin{center}
  \setlength{\unitlength}{1mm}
  \begin{picture}(50,30)(5,5)
   \put(10,10){\line(1,0){10}}
   \put(10,10){\line(0,1){10}}
   \put(20,10){\line(3,2){25}}
   \put(10,10){\circle*{1}}
   \put(10,20){\circle*{1}}
   \put(20,10){\circle*{1}}
   \put(40,23.5){\circle*{1}}
   \put(20,20){\circle*{1}}
   \put(-1,2){\makebox(10,10){$(0,0)$}}
   \put(17,2){\makebox(20,10){$w_1=(1,0)$}}
   \put(-1,15){\makebox(10,10){$(0,1)$}}
   \put(31,14){\makebox(30,10){$(a\beta,a)$}}
   \put(14,18){\makebox(10,10){$(1,1)$}}
    \put(25,24){\makebox(10,10){$G$}}
    \end{picture}
   \end{center}
   \caption{$G$ near singular $w_1$}
   \label{fig2g}
   \end{figure}
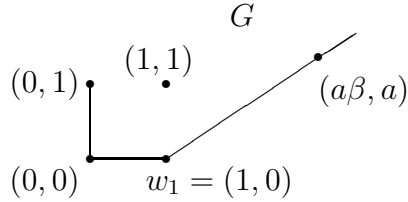

Finally, consider the case that $w_0$ is a nonsingular vertex and the next vertex $w_1$
 of $G$ is singular.  As before, set $w_0$ the origin, $w_1$ lying on the $x$-axis and
 the other edge from $w_0$ lying on the $y$-axis.
 Of course, the point $(1,1)$ is a vertex of $G^{\circ}$.
 If the $x$-coordinate of $w_1$ is greater than one, then the point $(2,1)$ is contained in 
 $G^{\circ}$ since $\dim G^{\circ}=2$.
 Set $w_1=(1,0)$.  See the Figure~\ref{fig2g}.
 
 %\noindent
 If $w_1$ be a  nonsingular vertex of $G$, then the other edge would have the direction
 $(\beta,1)$.  Since we assume that $w_1$ is  a singular vertex of type $(a,1)$
 with $a\ge2$ as in the Figure~\ref{fig2e} (b), 
  the edge from $w_1$ has the direction
 $(a\beta -1, a)$ because the direction is $a(\beta,1) +(-1,0)$ with respect to
 the basis $\{(-1,0), (\beta, 1)\}$.
 We see $\beta\ge1$ since $G^{\circ}$ contains the point $(1,1)$.
 If $\beta\ge2$, then the point $(2,1)$ is contained in $G^{\circ}$ since $a\beta -1>a$.
 Set $\beta=1$.  We distinguish two cases:
 
 \noindent
 If $a\ge3$, then $G^{\circ}$ contains the point $(2,2)$ as its boundary point since
 lattice points in $G$ are exhausted by the lines $y=x+k$ with $k\ge -1$.
 
 \noindent
 Set $a=2$.  If the point $(2,2)\in \partial G$ is a nonsingular vertex of $G$, then  the other edge from 
 $(2,2)$ does not have the direction $(0,1)$ since $\dim G^{\circ}=2$. 
 Thus $G^{\circ}$ contains the point $(2,3)$ as its boundary point by
 the exhausting the lines $y=2x+k$ with $k\ge -2$.
 
 \noindent
 If $w_r$ is a vertex of $G$, we see that the point $(1,2), (2,2)$ or $(3,2)$
 is contained in $G^{\circ}$ in the same way as $w_1$.
 Since the vertex $(1,1)$ of $G^{\circ}$ connects with two of the points
 $(1,2), (2,2), (2,3)$ and $(3,2)$, the vertex $(1,1)$ is nonsingular.
\hfill $\Box$

From this Lemma we see that facets of $Q$ have at worst $A_1$-singularities.

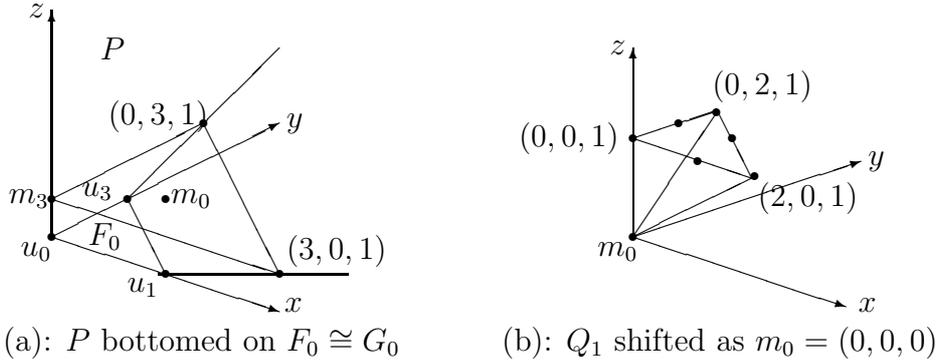
\begin{figure}[h]
 \begin{center}
 \begin{tabular}{lr}
  \setlength{\unitlength}{1mm}
  \begin{picture}(60,40)(5,0)
   \put(10,10){\vector(3,-1){30}}
   \put(10,10){\vector(2,1){30}}
   \put(10,10){\vector(0,1){30}}
   \put(20,15){\line(1,-2){5}}
   \put(20,15){\line(1,1){20}}
   \put(24,5){\line(1,0){25}}
   \put(30,25){\line(1,-2){10}}
   \put(10,15){\line(2,1){20}}
   \put(10,15){\line(3,-1){30}}
   \put(25,15){\circle*{1}}
   \put(21,10){\makebox(10,10)[r]{$m_0$}}
   \put(12,5){\makebox(10,10){$F_0$}}
   \put(3,3){\makebox(10,10){$u_0$}}
   \put(2,10){\makebox(10,10){$m_3$}}
   \put(17,2){\makebox(10,10)[b]{$u_1$}}
   \put(34,6){\makebox(20,10)[br]{$(3,0,1)$}}
   \put(14,21){\makebox(20,10){$(0,3,1)$}}
   \put(11,11){\makebox(10,10){$u_3$}}
   \put(13,30){\makebox(10,10){$P$}}
   \put(10,10){\circle*{1}}
   \put(25,5){\circle*{1}}
   \put(10,15){\circle*{1}}
   \put(20,15){\circle*{1}}
   \put(40,5){\circle*{1}}
   \put(30,25){\circle*{1}}
   \put(33,0){\makebox(10,10)[br]{$x$}}
   \put(3,35){\makebox(10,10){$z$}}
   \put(37,20){\makebox(10,10){$y$}}
  \end{picture}&
  
  \setlength{\unitlength}{1mm}
  \begin{picture}(50,40)(0,0)
  \put(10,10){\vector(3,-1){28}}
  \put(10,10){\vector(3,1){30}}
  \put(10,10){\vector(0,1){25}}
  \put(10,10){\line(2,1){16}}
  \put(10,10){\line(2,3){11}}
  \put(10,23){\line(3,-1){16}}
  \put(10,23){\line(3,1){11}}
  \put(21,27){\line(1,-2){4.5}}
  \put(10,10){\circle*{1}}
  \put(10,23){\circle*{1}}
  \put(21,26.5){\circle*{1}}
  \put(26,18){\circle*{1}}
  \put(16,25){\circle*{1}}
  \put(18.5,20){\circle*{1}}
  \put(23,23){\circle*{1}}
  \put(3,3){\makebox(10,10){$m_0$}}
  \put(-5,18){\makebox(20,10)[l]{$(0,0,1)$}}
  \put(23,13){\makebox(20,10)[b]{$(2,0,1)$}}
  \put(17,25){\makebox(20,10){$(0,2,1)$}}
  \put(32,0){\makebox(10,10)[br]{$x$}}
  \put(37,15){\makebox(10,10){$y$}}
  \put(3,30){\makebox(10,10){$z$}}
\end{picture}\\
\mbox{(a): $P$ bottomed on $F_0\cong G_0$} & \mbox{(b): $Q_1$ shifted as $m_0=(0,0,0)$}
\end{tabular}
 \end{center}
\caption{$P$ containing $Q_1$}
\label{fig4}
\end{figure}

Consider an example that $Q$ has a singularity at $m_0$.  
Assume that $P$ is locally described as $\{x\ge0, y\ge0, z\ge0, 2z\ge x+y-1\}$.  See Figure~\ref{fig4} (a).
Then $Q$ has a vertex $m_0=(1,1,1)$ and three edges of directions $(0,0,1), (2,0,1), (0,2,1)$.
See the Figure~\ref{fig4} (b).
Moreover we see that $Q$ has three facets meeting at $(1,1,1)$ as a singular vertex.
Thus after a suitable affine transformation of $M$,
we see that this $Q$ has the shape like
\begin{equation}\label{sect3:eq1}
Q_1:=\mbox{Conv}\{(0,0,0), (2,0,1), (0,2,1), (0,0,1)\}
\end{equation}
with $m_0=(0,0,0)$.
We call  $Q$ has a singularity of type $Q_1$ at $m_0$ in this case.

We note that $Q$ has three facets meeting at $m_0$ with $A_1$-singularity, that is, all two
of the three primitive points on the three edges have one lattice point between them.
We see that they are $(1,0,1), (0,1,1)$ and $(1,1,1)$ in the Figure~\ref{fig4} (b).

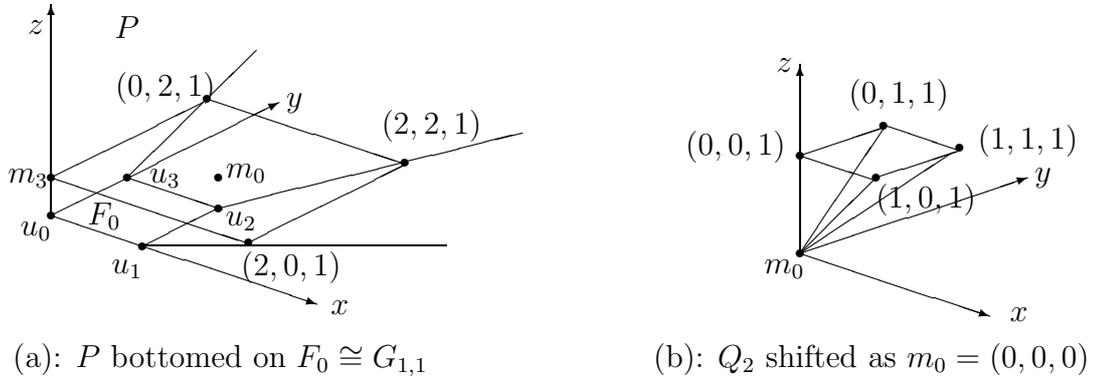
\begin{figure}[h]
 \begin{center}
  \setlength{\unitlength}{1mm}
  \begin{tabular}{lr}
  \begin{picture}(80,45)(5,-5)
   \put(10,10){\vector(3,-1){35}}
   \put(10,10){\vector(2,1){30}}
   \put(10,10){\vector(0,1){28}}
   \put(10,15){\line(3,-1){26}}
   \put(10,15){\line(2,1){20}}
   \put(22,6){\line(2,1){10}}
   \put(20,15){\line(3,-1){12.5}}
   \put(22,6){\line(1,0){40}}
   \put(20,15){\line(1,1){17}}
   \put(35,6){\line(2,1){22}}
   \put(30.5,25.5){\line(3,-1){25.5}}
   \put(32,11){\line(4,1){40}}
   \put(10,10){\circle*{1}}
   \put(10,15){\circle*{1}}
   \put(20,15){\circle*{1}}
   \put(22,6){\circle*{1}}
   \put(30.5,25.5){\circle*{1}}
   \put(32,11){\circle*{1}}
   \put(36,6.5){\circle*{1}}
   \put(56.5,17){\circle*{1}}
   \put(3,3){\makebox(10,10){$u_0$}}
   \put(15,-2){\makebox(10,10){$u_1$}}
   \put(43,-3){\makebox(10,10)[b]{$x$}}
   \put(3,30){\makebox(10,10){$z$}}
   \put(37,20){\makebox(10,10){$y$}}
   \put(2,10){\makebox(10,10){$m_3$}}
   \put(35,1){\makebox(20,10)[bl]{$(2,0,1)$}}
   \put(50,17){\makebox(20,10){$(2,2,1)$}}
   \put(15,22){\makebox(20,10){$(0,2,1)$}}
   \put(12,5){\makebox(10,10){$F_0$}}
   \put(20,10){\makebox(10,10){$u_3$}}
   \put(27,8){\makebox(10,10)[br]{$u_2$}}
   \put(15,30){\makebox(10,10){$P$}}
   \put(32,15){\circle*{1}}
   \put(33,15){\makebox{$m_0$}}
  \end{picture}&
  
  \begin{picture}(40,30)(8,0)
  \put(10,10){\vector(3,-1){25}}
  \put(10,10){\vector(3,1){30}}
  \put(10,10){\vector(0,1){25}}
   \put(10,10){\line(1,1){10}}
  \put(10,10){\line(2,3){11}}
  \put(10,10){\line(3,2){20}}
  \put(10,23){\line(3,-1){10}}
  \put(10,23){\line(3,1){11}}
  \put(21,27){\line(3,-1){10}}
  \put(20,20){\line(3,1){11}}
  \put(10,10){\circle*{1}}
  \put(10,23){\circle*{1}}
  \put(21,27){\circle*{1}}
  \put(20,20){\circle*{1}}
  \put(31,24){\circle*{1}}
  \put(3,3){\makebox(10,10){$m_0$}}
  \put(30,1){\makebox(10,10)[br]{$x$}}
  \put(37,15){\makebox(10,10){$y$}}
  \put(3,30){\makebox(10,10){$z$}}
  \put(-5,22){\makebox(20,10)[bl]{$(0,0,1)$}}
  \put(13,26){\makebox(20,10){$(0,1,1)$}}
  \put(30,20){\makebox(20,10){$(1,1,1)$}}
 %\put(15,0){\makebox(10,10)[b]{($Q_2$)}}
  \put(20,15){\makebox(20,10)[bl]{$(1,0,1)$}}
  \end{picture}\\
  \mbox{(a): $P$ bottomed on $F_0\cong G_{1,1}$} & \mbox{(b): $Q_2$ shifted as $m_0=(0,0,0)$}
\end{tabular}
 \end{center}
\caption{$P$ containing $Q_2$ }
\label{fig6}
\end{figure}

We have another example whose singular vertex $m_0$ is not singular in proper facets.
Assume that $P$ is locally described as $\{0\le x\le z+1, 0\le y\le z+1, z\ge0\}$. 
 See Figure~\ref{fig6} (a).
Then $Q$ has a vertex $m_0=(1,1,1)$ and four edges meeting at $m_0$ 
with directions $(0,0,1), (1,0,1),
(0,1,1), (1,1,1)$.  See the Figure~\ref{fig6} (b).
We see that  four facets of $Q$ meeting at $m_0$ are nonsingular at $m_0$.
After a suitable affine transformation we may draw the shape of $Q$ at $m_0$ like
\begin{equation}\label{sect3:eq2}
Q_2:= \mbox{Conv}\{(0,0,0), (1,0,1), (0,1,1), (0,0,1), (1,1,1)\}
\end{equation}
with $m_0=(0,0,0)$.
We call $Q$ has a singularity of type $Q_2$ at $m_0$ in this case.

We recall how to define a point $\bar{l}_{\sigma}$ in $P\cap M$   in this section.
For a vertex $u_0$ of $P$, we choose the cone $\sigma\in \Delta(3)$ so that
$\mathbb{R}_{\ge0}(P-u_0) = \sigma^{\vee}$.
If $\{m_1-u_0, m_2-u_0, m_3-u_0\}$ is the generator of the semi-group $\sigma^{\vee}
\cap M$, then we set $\bar{l}_{\sigma} = u_0 +\sum_{i=1}^3 (m_i-u_0)$.
Since this point  $\bar{l}_{\sigma}$ is defined only by the vertex $u_0$ of $P$,
we may write as $\bar{u}_0$.

It may happens $\bar{u}_0 =\bar{u}_1$ for different vertices $u_0$ and $u_1$ of $P$.
If $\bar{u}_i \not= \bar{u}_j$ for $i\not=j$, then $Q$ is nonsingular since every facet of $Q$
corresponds to a parallel facet of $P$ and when three facets of $Q$ meet at $\bar{u}$ 
corresponding three parallel facets of $P$ meet at $u$.

Before describing the singularity of $Q$, we will classify the facet $F_0$ of $P$
such that it happens $\bar{u}_0=\bar{u}_1$ for vertices $u_0\not=u_1$ of $F_0$.
As before, we set $u_0=(0,0,0)$ a vertex of $F_0$, $m_1=(1,0,0), m_2=(0,1,0)$
lying on the two edges of $F_0$ meeting at $u_0$ and $m_3=(0,0,1)$ lying on the other edge
of $P$ connecting with $u_0$.
Then $\bar{u}_0=(1,1,1)$.

\begin{figure}[h]
 \begin{center}
   \setlength{\unitlength}{1mm}
  \begin{picture}(60,30)(0,2)
   \put(10,10){\vector(3,-1){40}}
   \put(10,10){\vector(2,1){30}}
   \put(10,10){\vector(0,1){20}}
   \put(25,5){\line(1,0){15}}
   \put(25,5){\line(1,1){5}}
   \put(25,17.5){\circle*{1}}
   %\put(21,10){\makebox(10,10)[r]{$m_0$}}
   \put(15,5){\makebox(10,10){$F_0$}}
   \put(-8,3){\makebox(10,10){$u_0=(0,0,0)$}}
   %\put(2,10){\makebox(10,10){$m_3$}}
   \put(17,-2){\makebox(10,10)[b]{$u_1=(a,0,0)$}}
   \put(42,2){\makebox(20,10)[br]{$(a+c,0,1)$}}
   \put(30,7){\makebox(20,10){$(a+b,1,0)$}}
   \put(22,15){\makebox(10,10)[l]{$u_r=(0,a',0)$}}
   \put(13,20){\makebox(10,10){$P$}}
   \put(10,10){\circle*{1}}
   \put(25,5){\circle*{1}}
   \put(30,10){\circle*{1}}
   %\put(20,15){\circle*{1}}
 \put(40,5){\circle*{1}}
   %\put(30,25){\circle*{1}}
   \put(43,-5){\makebox(10,10)[br]{$x$}}
   \put(3,25){\makebox(10,10){$z$}}
   \put(37,20){\makebox(10,10){$y$}}
  \end{picture}
  \end{center}
  \caption{$P$ near the vertex $u_1$}
  \label{fig5a}
  \end{figure}
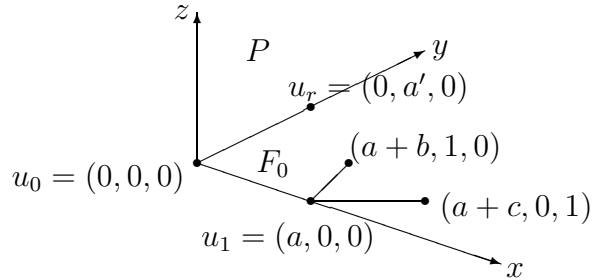

Let $u_1$ and $u_r$ be the vertices of $F_0$ adjacent to $u_0$.  Set $u_1=(a,0,0)$ and 
$u_r=(0,a',0)$.  From the vertex $u_1$, $P$ has the other two edges with the directions
$(b,1,0)$ and $(c,0,1)$ as in the Figure~\ref{fig5a}.  
We see $b, c\ge \max\{-1, -a+1\}$ since all facets of $P$ are
nonsingular.
Then $\bar{u}_1 = u_1+(-1,0,0)+(b,1,0)+(c,0,1)=(a+b+c-1,1,1)$.
Since $\bar{u}_0=(1,1,1)$,  we see that $\bar{u}_0=\bar{u}_1$ if and only if $a+b+c=2$.

\begin{lem}\label{sect3:l0}
In the above notation, if $\bar{u}_0=\bar{u}_1$, then we distinguish three cases (assuming
$b\le c$).
\begin{itemize}
\item[{\rm (a)}] $a=1$: When $b=-1$ and $c=2$, we see that $F_0$ is isomorphic to the basic
triangle $G_0=\mbox{\rm{Conv}}\{(0,0), (1,0), (0,1)\}$.
When $b=0$ and $c=1$, $F_0$ has two parallel edges of distance one.
\item[{\rm (b)}] $a=2$:  In this case, $b=-1$ and $c=1$.  $F_0$ is isomorphic to $2G_0$ or
the tetragon $G_{1,2}=\mbox{\rm{Conv}}\{(0,0), (0,1), (1,1), (2,0)\}$.
\item[{\rm (c)}] $a=3$: In this case, $b=-2$ and $c=1$.  $F_0$ is isomorphic to the teragon
$G_{1,3}=\mbox{\rm{Conv}}\{(0,0), (0,1), (1,1), (3,0)\}$.
\end{itemize}

Moreover, if $\bar{u}_r=\bar{u}_0$, then when $b=0$ in the case {\rm (a)} we see that
$F_0$ is isomorphic to the basic square $G_{1,1}=\mbox{\rm{Conv}}\{(0,0), (0,1), (1,1), (1,0)\}$.
In the other cases, we see $\bar{u}_r=\bar{u}_0$.

If $b>0$, then these $F_0$ appear as a facet of $P$ contained in the plane $\{y=0\}$.
\end{lem}

{\it Proof.}
When $a=1$, $b+c=1$ since $a+b+c=2$.  Since $c\ge b\ge -1$, we see $b=0$ or $b=-1$.

When $a=2$, $b+c=0$.  Since $c\ge b\ge -1$, $b=0$ or $b=-1$.  
If $b=0$, then $c=0$, hence $P$ is surrounded by the parallel two
planes $\{x=0\}$ and $\{ x=2\}$.  It does not happen since $\dim Q=3$.
Thus $b=-1$ and $c=1$.

When $a\ge3$, $b+c\le -1$.  Since $c\ge b$, $b<0$.
In this case, $F_0$ is contained in the triangle $T_0: 0\le x, 0\le y, x-by\le a$.
If $c\le0$, then $G$ is contained in $tT_0$ with $0<t\le1$ and 
$\mbox{Int} (tT_0)\cap M = \{(1,1,1)\}$.  This contradicts to $\dim Q=3$.
We see $c\ge1$.
Since $a+b\ge1, a+b+c=2$ and $c\ge1$, we see $a+b=1$ and $c=1$.
Then $b\le -2$ and $F_0\cong G_{1,a}$ since $F_0$ is nonsingular.

Consider the facet $F'$ of $P$ containing $u_1$ spanned by the directions $(b,1,0)$ and $(c,0,1)$.
On the facet $F'$ there is the lattice point $(a,0,0)+(b,1,0)+(c,0,1)=(2,1,1)$.  If the point $(2,1,1)$
is contained in the edge from $u_2=(a+b,1,0)\in F_0$, then $P$ cannot contain the point
$\bar{u}_0=(1,1,1)$ in its interior.
Hence, the lattice point  $(a,0,0)+2(b,1,0)+(c,0,1)=(a+2b+c,2,1)=(b+2,2,1)$ is on the facet $F'$.
Since this point is also contained in $P$, $b+2\ge0$.  This implies that $b=-2$, hence, $a=3$.
\hfill $\Box$

From this Lemma we see the singularity of $Q$.

\begin{prop}\label{sect3:p2}
Let $P$ be an integral convex polytope in $M_{\mathbb{R}}$ corresponding to a pair
$(X, L)$ of a nonsingular toric 3-fold $X$ and an ample line bundle $L$ on $X$
with $h^0(L\otimes \mathcal{O}_X(K_X))\not=0$.
Let $Q= \mbox{\textup{Conv}}(\mbox{\textup{Int}}(P)\cap M)$ be the interior polytope of
 $P$.   If $\dim\ Q=3$, then the singularities of $Q$ are the singular points of
the cones over $(\mathbb{P}^2, \mathcal{O}(2))$ and $(\mathbb{P}^1\times \mathbb{P}^1,
\mathcal{O}(1,1))$, which are given by the polytopes $Q_1$ in the Figure~\ref{fig4} {\rm(b)} and $Q_2$ 
in the Figure~\ref{fig6} {\rm(b)}, respectively.
\end{prop}

{\it Proof.}
Let $m_0$ be a vertex of $Q$.  Then we can choose a face $F_0$ of $P$ 
with the vertex $m_0$ such that the primitive elements $m_1, m_2\in F_0\cap M$ and
$m_3\in P\cap M$ on three edges meeting at $m_0$ form a $\mathbb{Z}$-basis of $M$
and that $m_0-u_0=(m_1-u_0)+(m_2-u_0)+(m_3-u_0)$ as in Figure~\ref{fig1}.
Let $(x,y,z)$ be the coordinates with respect to $\{(m_1-u_0), (m_2-u_0), (m_3-u_0)\}$.
Set $G=\{z=1\}\cap P$.   We may set $u_0=(0,0,0)$.

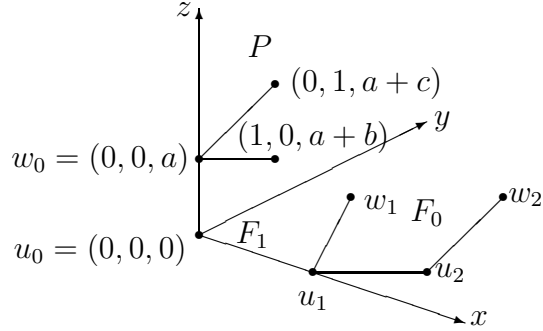
\begin{figure}[h]
 \begin{center}
   \setlength{\unitlength}{1mm}
  \begin{picture}(60,40)(0,5)
   \put(10,10){\vector(3,-1){35}}
   \put(10,10){\vector(2,1){30}}
   \put(10,10){\vector(0,1){30}}
   \put(10,20){\line(1,0){10}}
   \put(10,20){\line(1,1){10}}
   \put(25,5){\line(1,2){5}}
  \put(25,5){\line(1,0){15}}
  \put(40,5){\line(1,1){10}}
   \put(10,20){\circle*{1}}
   \put(20,18){\makebox(10,10){$(1,0,a+b)$}}
   \put(35,8){\makebox(10,10){$F_0$}}
   \put(12,5){\makebox(10,10){$F_1$}}
   \put(-8,3){\makebox(10,10){$u_0=(0,0,0)$}}
   \put(-8,15){\makebox(10,10){$w_0=(0,0,a)$}}
   \put(20,0){\makebox(10,10)[b]{$u_1$}}
   \put(22,28){\makebox(20,10)[br]{$(0,1,a+c)$}}
   \put(24,9){\makebox(20,10){$w_1$}}
   \put(38,0){\makebox(10,10){$u_2$}}
    \put(48,10){\makebox(10,10){$w_2$}}
   \put(13,30){\makebox(10,10){$P$}}
   \put(10,10){\circle*{1}}
   \put(25,5){\circle*{1}}
   \put(30,15){\circle*{1}}
   \put(20,20){\circle*{1}}
 \put(40,5){\circle*{1}}
   \put(20,30){\circle*{1}}
   \put(50,15){\circle*{1}}
   \put(38,-2){\makebox(10,10)[br]{$x$}}
   \put(3,35){\makebox(10,10){$z$}}
   \put(37,20){\makebox(10,10){$y$}}
  \end{picture}
  \end{center}
  \caption{$P$ near the vertex $w_0$}
  \label{fig5b}
  \end{figure}

From Lemma~\ref{sect3:l0}, we have a classification of the facets $F_0$ of $P$ such that 
$G^{\circ}$ is one point set $\{(1,1,1)\}$.
Set $F_1:=P\cap\{y=0\}$ and $F_2:=P\cap\{x=0\}$.
Let $w_0=(0,0,a)$ with $a\ge1$ be the vertex of $P$ connecting with $u_0$.
From the vertex $w_0$, $P$ has the other two edges with the directions $(1,0,b)$ and $(0,1,c)$.
See the Figure~\ref{fig5b}.

We distinguish three cases according to the cases in Lemma~\ref{sect3:l0}.

%\noindent
Case (a): $F_0\cong G_0$ or $G_{1,1}$.  In this case, $G\cong 3G_0$ or $2G_{1,1}$.
When $a=1$ we may assume $b, c\ge0$ since $Q\not=\emptyset$.
If $b+c\ge2$, then $\bar{w}_0\not=\bar{u}_0$.
If $b=0, c=1$ and $F_0\cong G_{1,1}$, then $P$ cannot contain $Q$ with $\dim Q=3$.
If $b=0, c=1$ and $F_0\cong G_0$, then $\bar{w}_0=\bar{u}_0$ and this is the case (c)
by exchanging the role of $F_0$ with the facet $F_2$.
When $a=2$ or $3$, we see $b\not=-1$ from the shape of $F_1$ or $F_2$.
From Lemma~\ref{sect3:l0}, $\bar{w}_0\not=\bar{u}_0$.
Except the case (c), the vertex $\bar{u}_0$ has three or four edges in $Q$ parallel to 
the edges of $P$ meeting at the vertices of $F_0$.
This $\bar{u}_0$ is a singular vertex of $Q$ of type $Q_1$ or $Q_2$.

%\noindent
Case (b):  $F_0\cong 2G_0$ or $G_{1,2}$.  In this case, $G\cong 3G_0$ or 
$\mbox{Conv}\{(0,0), (0,2), (1,2), (3,0)\}$.
When $1\le a\le3$ we may assume $b, c\ge0$ from the shape of nonsingular facets
$F_1$ and $F_2$.
When $a=1$ if $b=1$ or $c=1$, then $\dim Q=1$.
From Lemma~\ref{sect3:l0}, $\bar{w}_0\not=\bar{u}_0$.
Thus if $F_0\cong 2G_0$, then the vertex $\bar{u}_0$ has three edges in $Q$ 
parallel to the three edges of $P$ meeting the vertices of $F_0$,
which have the directions $(0,0,1), (1,0,1),(0,1,1)$.
In this case $(1,1,1)$ is a nonsingular vertex of $Q$.

Even in the case $F_0\cong G_{1,2}$, the interior polytope $Q$ has tow
edges from the vertex $\bar{u}_0$ with the directions $(0,0,1)$ and $(1,0,1)$.
Set $u_2=(1,1,0)$ and $u_3=(0,1,0)$ the vertices of $F_0$.
In this case we see $\bar{u}_2=\bar{u}_3$ from Lemma~\ref{sect3:l0} (a).
Set $F'=P\cap\{y=z+1\}$ the facet containing $u_2$ and $u_3$.
On the edges of $F'$ we take vertices $w_2=(1,d'+1,d'), w_3=(0,d+1,d)$.
If $d, d'\ge2$, then $\bar{w}_i\not=\bar{u}_i$ from Lemma~\ref{sect3:l0},
 and $Q$ has an edgs meeting $\bar{u}_0$ with the direction $(0,1,1)$.
If $d=d'=1$, then $P$ is surrounded by the parallel twi planes $\{y=0\}$ and $\{y=2\}$,
hence, $P$ cannot contain $Q$ with dimension three.
Even if $d=1$ and $d'\ge2$, then $P$ contains the point $(1,2,2)$ in its interior
because $Q$ already contains the point $(1,1,2)$.
Since $P$ is contained in the combinatorial prism $\{0\le x, 0\le y, 0\le z, x+y+2\le z\}$
and since the section of the prism at $\{z=2\}$ contains only three points $(1,1,2), (2,1,2),
(1,2,2)\}$ in its interior, the point $(1,1,1)$ is a nonsingular vertex of $Q$ spanned by
vectors $(0,0,1), (1,0,1),(0,1,1)$.
If $d\ge2$ and $d'=1$, we obtain the same conclusion by exchanging the role of $u_0, u_1$
and $u_2, u_3$.

%\noindent
Case (c): $F_0\cong G_{1,3}$.  In this case, $G \cong\mbox{Conv}\{(0,0), (0,2),(4,0)\}$
and $(0,2,1)$ is a vertex of $P$.
Moreover, $P$ has the facet $\mbox{Conv}\{(0,1,0), (1,1,0),(0,2,1)\}\cong G_0$
and the third edge of $P$ from $(0,2,1)$ passes through $(0,3,3)$.
As before set $w_0=(a,0,0)$ the vertex of $P$.  When $1\le a\le3$, we see $b, c\ge0$
from the shape of facets $F_1$ and $F_1$.  If $a=1, b=0$ and $c=1$, then
$F_1$ and $F_2$ are nonsingular, but $\dim Q=0$.
From Lemma~\ref{sect3:l0}, we see $\bar{w}_0\not=\bar{u}_0$.
For a vertex $w_1=(a+3,0,a)$ on the edge of $P$ meeting $u_1=(3,0,0)$, we see
$\bar{w}_1\not=\bar{u}_1$ in the same way.  Thus $Q$ has two edges from $(1,1,1)$
with the directions $(0,0,1)$ and $(1,0,1)$.

At the vertex $w_3:=(0,2,1)$ we have a basis $\{(0,-1,-1),(1,-1,-1),(0,1,2)\}$ 
from the primitive elements on
the three edges of $P$ meeting at the vertex.
 As in the case $w_0$, we consider the point $w=(0,2,1)+a(0,1,2)$ and two vectors $b(0,1,2)+(0,-1,-1)$
 and $c(0,1,2)+(1,-1,-1)$ from $w$ with respect to this basis.
 If $a=1$, then $b,c\ge1$ from the shape of the facet $F_1$.  If $a=2$, then two vectors $(0,-1,-1)$
 and $(1,-1,-1)$ define the plane $\{z=y+1\}$ (the case $b=c=0$), hence, we see $\dim Q=0$.
 If $a=3$, then $b,c\not=-1$ from the shape of $F_1$.
From Lemma~\ref{sect3:l0}, we see $\bar{w}\not=\bar{u}_0=\bar{w}_3$.
Then $Q$ has the edge from $(1,1,1)$ with the direction $(0,1,2)$, hence, the vertex $\bar{u}_0=
(1,1,1)$ is nonsingular.
\hfill $\Box$

\section{Adjacent singular vertices.}\label{sect5}

In the previous section we investigate the singularities of $Q$ when $\dim Q=3$.  
They are two types of singularities
described by the polytopes $Q_1$ and $Q_2$.
In this section we treat the case that $Q$ has an edge whose ends are the singularities of $Q$
and which contains no more lattice points.

Consider the case that $Q$ has the singularity of  type  $Q_1$ at $m_0$.
In this case $P$ is locally of the form as in the Figure~\ref{fig4} (a). 
In the Figures~\ref{fig4} (b) and \ref{fig6} (b), 
we see that the part $Q\cap\{0\le z\le1\}$ of $Q$ containing the singular vertex
$m_0$ is normally generated by decomposing into a union of basic 3-simplices.
We need a way to decompose $Q$ into a union of normally generated polytopes under some condition.

\begin{lem}\label{sect5:l1}
Let $m_0=(0,0,0)$ be the singular vertex of $Q$.
Assume that the point $(0,0,1)$ in the Figure~\ref{fig4} {\rm(b)} or \ref{fig6} {\rm(b)} is a vertex of $Q$.
If  $(1,0,a)$ and $(0,1,b)$ $(1\le a\le b)$ are contained in $Q$, then
$a=1$ and $b\le3$, or $a=b=2$.

Moreover, if $(0,0,1)$ is a singular vertex of type $Q_1$, then $Q$ has the other two
edges meeting at $(0,0,1)$ connecting $(2,0,2)$ and $(0,2,2)$, respectively.
If $(0,0,1)$ is a singular vertex of type $Q_2$, then $Q$ has the other three edges
meeting at $(0,0,1)$ connecting $(1,0,1), (0,1,1)$ and $(1,1,2)$, respectively.
\end{lem}

{\it Proof.}
See the Figures~\ref{fig4} (a) and \ref{fig6} (a). 
Set $w_0=(0,0,d)$ the vertex of $P$ and $(1,0,f), (0,1,g)$ the directions of two
other edges of $P$ from $w_0$.  We may assume $f\le g$.
Set $F_1=P\cap\{y=0\}$ and $F_2=P\cap\{x=0\}$ the facets of $P$.

We note that $P$ contains the points $(2,1,2), (1,2,2)$ in its interior by assumption.
If $\bar{w}_0=(1,1,2)$, then $d+f+g=3$ since $\bar{w}_0 = w_0 +(1,0,f)+(0,1,g)
+(0,0,-1)$.
In this case we note $f\ge0$ since $2(1,0,f)+(0,1,g)+w_0=(2,1,3+f)$.

\noindent
When $d=3$, we have $f=g=0$ and $Q$ is contained in $\{0\le z\le1\}$.

\noindent
When $d=2$,  we have $f=0, g=1$.  In this case, since $P$ is contained in
$\{z\le y+2\}$ we have $a=1, b\le2$.

\noindent
When $d=1$, we have two cases $f=0,g=2$ and $f=g=1$.
If $f=0$ and $g=2$, then we see $a=1$ and $b\le3$ since $P$ is contained in
$\{z\le 2y+1\}$.
If $f=g=1$, then we see $a, b\le2$ since $P$ is contained in $\{z\le x+y+1\}$.

Next we assume $\bar{w}_0=(1,1,2)$ is a singular vertex of $Q$  of type $Q_1$.
From Lemma~\ref{sect3:l0} and Proposition~\ref{sect3:p2},
the vertex $w_0=(0,0,d)$ is also a vertex of a facet $F'\cong G_0$ of $P$.
If $d=2$ and $f=0$, then the point $(1,0,2)$ has to be a vertex with the edge
of the direction $(2,0,-1)$, which connects to $(3,0,1)$.  
This is impossible since $F_1$ is nonsingular.
Thus we have $d=f=g=1$.
In this case, the vertex $(1,0,2)$ has the other edge with the direction $(2,0,1)$
and the vertex $(0,1,2)$ has the edge with the direction $(0,2,1)$.
From Lemma~\ref{sect3:l0} we see that the singular vertex $\bar{w}_0$ of $Q$
has two other edges with the directions $(2,0,1)$ and $(0,2,1)$.

Finally we assume that $\bar{w}_0$ is a singular vertex of type $Q_2$.
In the same way we see $d=f=g=1$.  From Lemma~\ref{sect3:l0} and 
Proposition~\ref{sect3:p2}, we see that the point $(1,1,3)$ is a vertex of $P$,
it has the other edge with the direction $(1,1,1)$,
the vertex $(1,0,2)$ has  the edge with the direction $(1,0,0)$ and
the vertex $(0,1,2)$ has the edge with the direction $(0,1,0)$.
From Lemma~\ref{sect3:l0}, we see the singular vertex $\bar{w}_0$ has the other three
edges with the directions $(1,0,1),(0,1,0)$ and $(1,1,1)$.
\hfill $\Box$

First we consider the case that $Q$ has an edge whose both ends are  singularities.

\begin{prop}\label{sect5:p11}
If $Q$ has an edge whose both ends are the singularities of $Q$ and which contains no
more lattice points, then $Q$ is  normally generated
unless the pair of the singularities consists of two $Q_1$'s.
\end{prop}

{\it Proof.}
(a) First we consider the case that $Q$ has the singularity of type $Q_2$ at the vertex  $m_0$.
 If $(0,0,1)$ in the Figure~\ref{fig6} (b)  is the singularity of $Q$ of type $Q_2$, then 
 the other three edges meeting at $(0,0,1)$ connect $(1,0,1), (0,1,1)$ and $(1,1,2)$
 from Lemma~\ref{sect5:l1}.

 In this case, the lattice points $(1,0,1)$ and $(0,1,1)$ are also vertices of $Q$.
 The other edges from the vertices $(1,0,1), (0,1,1)$ have the same direction $(1,1,1)$.
Hence $Q$ is contained in the quadrangular prism $\{x\le z\le x+1, y\le z\le y+1\}$
with the bottom
$\mbox{Conv}\{(0,0,1),$ $(1,0,1), (1,1,1), (0,1,1)\}$ and the direction vector $(1,1,1)$.

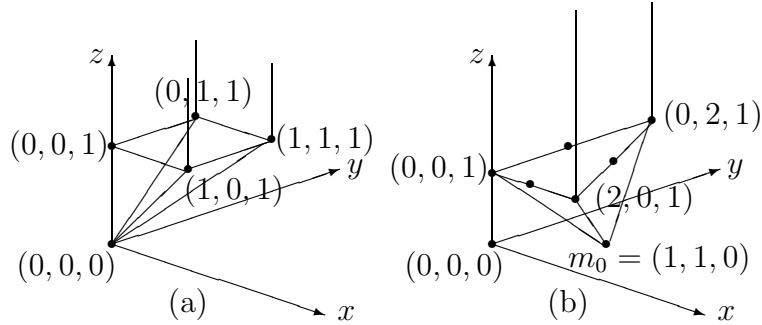
\begin{figure}[h]
 \begin{center}
  \setlength{\unitlength}{1mm}
  \begin{picture}(100,40)(0,0)
  \put(10,10){\vector(3,-1){28}}
  \put(10,10){\vector(3,1){30}}
  \put(10,10){\vector(0,1){25}}
  \put(10,10){\line(1,1){10}}
  \put(10,10){\line(2,3){11}}
  \put(10,10){\line(3,2){20}}
  \put(10,23){\line(3,-1){10}}
  \put(10,23){\line(3,1){11}}
  \put(21,27){\line(3,-1){10}}
  \put(20,20){\line(3,1){11}}
  \put(20,20){\line(0,1){12}}
  \put(21,27){\line(0,1){10}}
  \put(31,24){\line(0,1){10}}
  \put(10,10){\circle*{1}}
  \put(10,23){\circle*{1}}
  \put(21,27){\circle*{1}}
  \put(20,20){\circle*{1}}
  \put(31,24){\circle*{1}}
  \put(-6,2){\makebox(20,10){$(0,0,0)$}}
  \put(-7,18){\makebox(20,10){$(0,0,1)$}}
  \put(28,16){\makebox(20,10)[t]{$(1,1,1)$}}
  \put(12,25){\makebox(20,10){$(0,1,1)$}}
  \put(16,12){\makebox(20,10){$(1,0,1)$}}
  \put(32,0){\makebox(10,10)[br]{$x$}}
  \put(37,15){\makebox(10,10){$y$}}
  \put(3,30){\makebox(10,10){$z$}}
  %\put(10,12){\makebox(10,10){$Q_1$}}
  \put(15,0){\makebox(10,10)[b]{(a)}}
  \put(60,10){\vector(3,-1){28}}
  \put(60,10){\vector(3,1){30}}
  \put(60,10){\vector(0,1){25}}
  \put(71,16){\line(1,1){10}}
  \put(71,16){\line(0,1){25}}
  \put(60,19.7){\line(3,-2){15}}
  \put(60,19.5){\line(3,1){21}}
  \put(60,19.5){\line(3,-1){11}}
  \put(81,27){\line(-1,-3){5.5}}
  \put(81,27){\line(0,1){15}}
  \put(71,16){\line(2,-3){4}}
  \put(60,10){\circle*{1}}
  \put(60,19.5){\circle*{1}}
  \put(75,10){\circle*{1}}
  \put(81,26.5){\circle*{1}}
  \put(71,16){\circle*{1}}
  \put(70,23){\circle*{1}}
   \put(76,21){\circle*{1}}
   \put(65,18){\circle*{1}}
  \put(45,2){\makebox(20,10){$(0,0,0)$}}
  \put(82,0){\makebox(10,10)[br]{$x$}}
  \put(87,15){\makebox(10,10){$y$}}
  \put(53,30){\makebox(10,10){$z$}}
  \put(43,15){\makebox(20,10){$(0,0,1)$}}
  \put(79,22){\makebox(20,10){$(0,2,1)$}}
  \put(70,11){\makebox(20,10){$(2,0,1)$}}
 %\put(75,17){\makebox(20,10){$(2,2,2)$}}
 \put(65,0){\makebox(10,10)[b]{(b)}}
 \put(72,3){\makebox(20,10){$m_0=(1,1,0)$}}
  \end{picture}
 \end{center}
\caption{$Q$ with a pair of $Q_2$'s and $Q_1$'s}
\label{fig8}
\end{figure}

By a suitable affine transformation of $M$ 
we may write $Q$ in the region $\{0\le x\le1, 0\le y\le1, 0\le z\}$
as in the Figure~\ref{fig8} (a).
The opposite side of this prism has four vertices.  If all four vertices are on a plane, then 
they are nonsingular vertices.  If the opposite side has two facets, then it has
 two $Q_2$'s as singularity as in the original side since singularities of $Q$ are only of type $Q_1$ or $Q_2$.
 Thus we can cut off $Q_2$'s
so that the rest $Q'$ is nonsingular.
From Proposition~\ref{sect2:p2} and Lemma~\ref{sect2:l2} it is normally generated.

(b) Consider the case that $Q$ has the singularity of type $Q_1$ at $m_0$ as in the Figure~\ref{fig4} (b).
If $(0,0,1)$ in the Figure  is the singularity of $Q$ of type $Q_2$, then
the other three edges meeting at $(0,0,1)$ are connected with $(1,0,1), (0,1,1)$ and $(1,1,1)$
as in the case (a).  Then the edge through $(1,0,1)$ goes to $(2,0,1)$ and the edge through $(0,1,1)$
also goes to $(0,2,1)$.
These are both vertices of $Q$.  
Since $(2,0,1)$ is a nonsingular vertex of the facet $\mbox{Conv}\{(0,0,0), (0,0,1),(2,0,1)$
of $Q$, the vertex $(2,0,1)$ of $Q$ is nonsingular or singular of type $Q_2$.

If $(2,0,1)$ is the singularity of type $Q_2$, then the other two edges from $(2,0,1)$
have the directions $(-1,1,0)$ and $(0,1,1)$ from Lemma~\ref{sect5:l1} after a suitable
transformation of $M$.
Since one edge with the direction $(-1,1,0)$ connects with the point $(0,2,1)$, 
the vertex $(0,2,1)$ is also a singular vertex of type $Q_2$.
The other edge goes to $(2,1,2)$, and if $(2,1,2)$ is not a vertex, then the edge goes to $(2,2,3)$.
Thus the other edge from $(0,2,1)$ goes to $(1,2,2)$ and may go to $(2,2,3)$.

If the point $(2,1,2)$ is a vertex, then the point $(1,2,2)$ is also a vertex and they make
a facet with the vertex $(1,1,2)$ isomorphic to the basic triangle $G_0$.

If we decompose $Q$ into two parts by cutting
at the plane $\{z=1\}$, then the upper polytope is isomorphic to
the twice of the basic 3-simplex $2P_0$ or $P_1$, both are nonsingular polytopes without
lattice points in its interiors.
 In this case, we see that $Q$  is normally generated.
 
 If the point $(2,0,1)$ is a nonsingular vertex, then the point $(0,2,1)$ is also a nonsingular vertex from above
  and the other edge from $(2,0,1)$
 has the direction $(1,1,1)$.  The vertex $(0,2,1)$ also has the edge with the direction $(1,1,1)$.
Hence,  $Q$ is contained in the triangular prism $\{x+1\ge z, y+1\ge z, x+y\le z+1\}$
with the bottom
 $\mbox{Conv}\{(0,0,1),$ $(0,2,1), (2,0,1)\}$ and the direction vector $(1,1,1)$.
By exchanging the vector $(0,0,1)$ with $(1,1,1)$ as a part of the basis of $M$
and shifting with the direction $(1,1,0)$, 
we may write $Q$ in the region $\{x, y, z\ge0, x+y\le2\}$
as in the Figure~\ref{fig8} (b), where $m_0$  moves to $(1,1,0)$.
In the Figure, if we cut off the polytope $Q\cap\{z\le1\}\cong Q_1$, then
the bottom of the rest is nonsingular.

Consider the opposite side of this prism.  Set $w_0=(0,0,a), w_1=(2,0,b)$ and $w_2=(0,2,c)$
the vertices of the prism in the opposite side.
If all three line segments $\overline{w_0w_1}, \overline{w_1w_2}$ and $\overline{w_2w_0}$
are also the edges of $Q$, then we claim that there is one more lattice point on each edge.

\noindent
If all three edges contains only end points as their lattice points, then all differences $a-b, b-c$ and $c-a$
of the $z$-coordinates have to be odd, but it is impossible.
Then we see that one edge has one more lattice point on it.  Say it is $\overline{w_0w_1}$.
In this case, two ends $w_0$ and $w_1$ of the edge are nonsingular vertices of the facet
$Q\cap\{y=0\}$.  
Since $w_0$ and $w_1$ are simplicial vertices of $Q$, they are also nonsingular vertices of $Q$
from Proposition~\ref{sect3:p2}.  Then the other two edges have one more lattice point.

From this consideration, we see that if the opposite side of the prism consists of one facet,
then all three vertices of the facet are nonsingular.

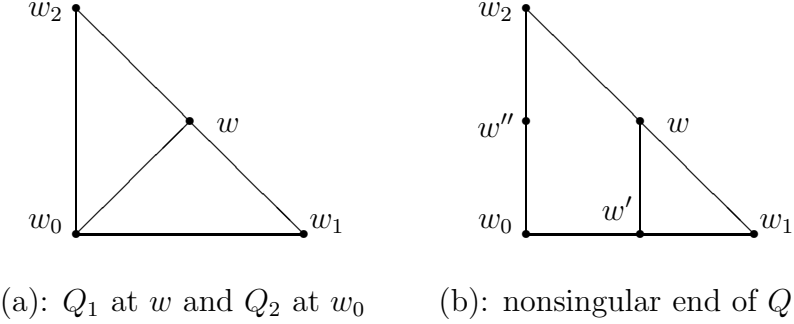
\begin{figure}[h]
\begin{center}
 \setlength{\unitlength}{1mm}
 \begin{tabular}{lr}
 \begin{picture}(50,45)(0,5)
  \put(10,10){\line(1,0){30}}
  \put(10,10){\line(0,1){30}}
  \put(10,10){\line(1,1){15}}
  \put(10,40){\line(1,-1){30}}
  \put(10,10){\circle*{1}}
  \put(10,40){\circle*{1}}
  \put(25,25){\circle*{1}}
  \put(40,10){\circle*{1}}
  \put(1,2){\makebox(10,19){$w_0$}}
  \put(38,2){\makebox(10,19){$w_1$}}
  \put(1,30){\makebox(10,19){$w_2$}}
  \put(25,15){\makebox(10,19){$w$}}
  \end{picture} &
  
  \begin{picture}(50,45)(-5,5)
   \put(10,10){\line(1,0){30}}
  \put(10,10){\line(0,1){30}}
  \put(25,10){\line(0,1){15}}
  \put(10,40){\line(1,-1){30}}
  \put(10,10){\circle*{1}}
  \put(10,25){\circle*{1}}
  \put(10,40){\circle*{1}}
  \put(25,25){\circle*{1}}
  \put(25,10){\circle*{1}}
  \put(40,10){\circle*{1}}
  \put(1,2){\makebox(10,19){$w_0$}}
  \put(38,2){\makebox(10,19){$w_1$}}
  \put(1,30){\makebox(10,19){$w_2$}}
  \put(25,15){\makebox(10,19){$w$}}
  \put(17,4){\makebox(10,19){$w'$}}
  \put(1,15){\makebox(10,19){$w''$}}
  \end{picture} \\
  \mbox{(a): $Q_1$ at $w$ and $Q_2$ at $w_0$} & \mbox{(b): nonsingular end of $Q$}
  \end{tabular}

\end{center}
\caption{opposite side of $Q$ with singularities $Q_1$ and $Q_2$}
\label{fig8e}
\end{figure}

If the point $w=(1,1,d)$ on the opposite side of $Q$ is a vertex connecting with $w_0$
(see the Figure~\ref{fig8e} (a)),
then $w_0$ is singular of type $Q_2$ (since it is not a simplicial vertex) and $w$ is also a
singular vertex of type $Q_1$.  By putting the prism upside down and taking a 
suitable transformation of $M$, we have the shape of $Q$ near $w$ as in the Figure~\ref{fig8} (b).

\noindent
To see this, set $w_0=(0,0,a), w_1=(2,0,a+2b)$ and $w_2=(0,2,a+2c)$.  Then
$w=(0,0,a-1)+(1,0,b+1)+(0,1,c+1)=(1,1,a+b+c+1)$ since $w$ is of type $Q_2$,
and we see that the line segment $\overline{w_1w_2}$ contains the lattice point $(1,1,a+b+c)$.
This implies $w$ is a singular vertex of type $Q_1$.

Moreover,
if the point $w'=(1,0,e)$ on the opposite side of $Q$ is a vertex connecting the vertex $w=(1,1,d)$,
then $w$ has to be a singular vertex of type $Q_2$, which is impossible since
$w$ is a singular vertex of the facet $\mbox{Conv}\{w_0,w, w_2\}$.

If the point $w'$  connects with $w$ and if $w_0$ is a nonsingular vertex
(see the Figure~\ref{fig8e} (b)), then $w'$ is a
nonsingular vertex of $Q$ from Proposition~\ref{sect3:p2} since  the edge $\overline{ww'}$ contains only two
lattice points.
Then $w$ is also nonsingular and so is $w_0$.

Moreover, if $w$ connects with the point $w''=(0,1,f)$, then $w$ has to be a singular vertex of type $Q_2$,
which is impossible.
To see this,  we may set as that vertices $w, w'$ and $w_1$ are on the plane $\{z=d\}$ and
$w_2=(0,2,d-1)$ since $w$ is a nonsingular vertex of the facet $Q\cap\{x+y=2\}$.
Since $w'$ is also a nonsingular vertex of the facet $Q\cap\{y=0\}$, we have $w_0=(0,0,d-1)$.
This implies that $w''$ is not vertex.

In any case, if we cut one or two polytopes isomorphic to  $Q_2$, then the rest is a nonsingular
polytope without lattice points in its interior.  Thus 
 we see that
 $Q$ is normally generated from Proposition~\ref{sect2:p2} and Lemma~\ref{sect2:l2}.
\hfill $\Box$

\begin{figure}[h]
 \begin{center}
  \setlength{\unitlength}{1mm}
  \begin{picture}(70,50)(0,0)
  \put(10,10){\vector(3,-1){28}}
  \put(10,10){\vector(3,1){30}}
  \put(10,10){\vector(0,1){25}}
  \put(10,10){\line(4,1){17}}
  \put(10,10){\line(1,1){15}}
  \put(27,14.5){\line(0,1){10}}
  \put(25,25){\line(0,1){10}}
  \put(10,20){\line(4,1){17}}
  \put(10,20){\line(1,1){15}}
  \put(27,24.5){\line(1,1){15}}
  \put(27,14.5){\line(1,1){15}}
  \put(25,25){\line(4,1){17}}
  \put(25,35){\line(4,1){17}}
  \put(42,29){\line(0,1){10}}
  \put(10,10){\circle*{1}}
  \put(10,20){\circle*{1}}
  \put(27,14.5){\circle*{1}}
  \put(25,25){\circle*{1}}
  \put(27,24.5){\circle*{1}}
  \put(42,39){\circle*{1}}
  \put(42,29){\circle*{1}}
  \put(25,35){\circle*{1}}
  \put(3,3){\makebox(10,10){$m_0$}}
  \put(-5,15){\makebox(20,10)[l]{$(0,0,1)$}}
  \put(21,8){\makebox(20,10)[b]{$(2,0,1)$}}
  \put(11,21){\makebox(20,10)[l]{$(0,2,1)$}}
  \put(28,19){\makebox(20,10)[l]{$(2,0,2)$}}
  \put(18,33){\makebox(20,10)[l]{$(0,2,2)$}}
  \put(43,23){\makebox(20,10)[l]{$(2,2,2)$}}
  \put(43,35){\makebox(20,10)[l]{$(2,2,3)$}}
  \put(32,0){\makebox(10,10)[br]{$x$}}
  \put(37,15){\makebox(10,10){$y$}}
  \put(3,30){\makebox(10,10){$z$}}
\end{picture}
 \end{center}
\caption{$Q$ with singularities  $Q_1$}
\label{fig9}
\end{figure}
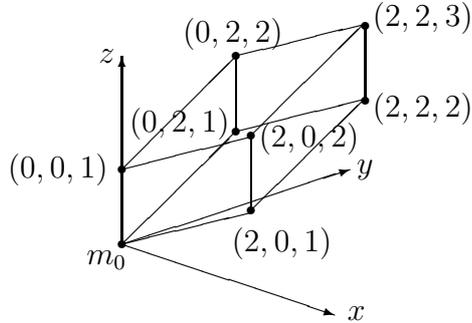

\begin{prop}\label{sect5:p2}  Let $Q$ be a singular polytope with the singularity at $m_0$
of type $Q_1$.
If all points $(0,0,1),(2,0,1),(0,2,1)$ of
 $Q$ in the Figure~\ref{fig4} are  singularities of type $Q_1$, 
 then $Q$ is normally generated.
\end{prop}
{\it Proof.}
In the Figure~\ref{fig4}, if $(0,0,1)$ is a singular vertex of type $Q_1$, then the other two edges
from $(0,0,1)$ connect with $(2,0,2)$ and $(0,2,2)$ from Lemma~\ref{sect5:l1}.
If the point $(2,0,1)$ is also a singular vertex of type $Q_1$, then then the other two edges
from $(2,0,1)$ connect with $(2,0,2)$ and $(2,2,2)$ from Lemma~\ref{sect5:l1}.
Then the point $(2,0,2)$ is also a singular vertex of type $Q_1$.
hence the polytope $Q$ satisfying the condition of the Proposition is a parallelotope with
eight vertices and all vertices are singular of type $Q_1$.  See the Figure~\ref{fig9}.

Then it is easy to see $Q$ is normally generated.
For example, cut $Q$ at the planes $\{z=1\}$ and $\{z=2\}$ into three pieces. 
The bottom and top parts are isomorphic to $Q_1$, hence, they are normally generated.
The rest $Q'$ has two parallel facets of the distance one.  If you cut it at the plane $\{x+y=2\}$, then
both have a nonsingular facet $R$ of a triangle isomorphic to $2G_0$ 
and a line segment $E$ parallel to one edge of the facet with distance one.
From Remark~\ref{sect2:rm1}, we see that both two pieces divided from  $Q'$ are normally generated
since we have
$$
R\cap \mathbb{Z}^2 + E\cap \mathbb{Z}^2 = (R+E)\cap \mathbb{Z}^2.
$$
This implies two teragonal cones are normally generated as in the proof of Lemma~\ref{sect2:l3}.
\hfill $\Box$

Next we need to consider the case that $Q$ has an edge whose one end is singular and
 another end is a nonsingular vertex.
Let  $m_0=(0,0,0)$ be a singular vertex of $Q$.
In  the Figures~\ref{fig4} (b) and \ref{fig6} (b),  we assume that $(0,0,1)$ is a nonsingular vertex of $Q$.
Lemma~\ref{sect5:l1} says that the other two edges from the vertex $(0,0,1)$ pass through
the points $(1,0,a)$ and $(0,1,b)$ with $1\le a\le3$ and $b=1$, or $a=b=2$.
Thus we may define the types of nonsingular vertices:
\begin{itemize}
 \item[(1,0)] The vertex $(0,0,1)$ is connected with $(1,0,2)$ and $(0,1,1)$.
 \item[(1,1)] The vertex $(0,0,1)$ is connected with $(1,0,2)$ and $(0,1,2)$.
  \item[(2,0)] The vertex $(0,0,1)$ is connected with $(1,0,3)$ and $(0,1,1)$.
\end{itemize}
For convenience of drawing pictures, we change to $a\ge b$ from the statement of Lemma~\ref{sect5:l1}.

By using this classification we have the following Proposition.

\begin{figure}[h]
 \begin{center}
   \setlength{\unitlength}{1mm}
 \begin{tabular}{lr}
  \begin{picture}(50,60)(-3,0)
  \put(10,10){\vector(3,-1){28}}
  \put(10,10){\vector(3,1){38}}
  \put(10,10){\vector(0,1){35}}
  \put(10,10){\line(1,1){20}}
  \put(10,10){\line(2,3){23}}
  \put(10,10){\line(3,2){20}}
  \put(10,23){\line(3,-1){10}}
  \put(10,23){\line(3,1){11}}
  \put(21,27){\line(3,-1){10}}
  \put(20,20){\line(3,1){11}}
  \put(20,33){\line(3,-1){10}}
  \put(20,33){\line(1,5){1.5}}
  \put(10,23){\line(1,1){10}}
  \put(10,23){\line(2,3){11}}
  \put(21,40){\line(3,1){12}}
  \put(31,24){\line(1,1){9}}
  \put(31,24){\line(2,3){10}}
  \put(30,30){\line(3,1){9.5}}
  \put(40,33){\line(1,5){1.2}}
  \put(41,39){\line(-3,2){8}}
  \put(10,10){\circle*{1}}
  \put(10,23){\circle*{1}}
  \put(21,27){\circle*{1}}
  \put(20,20){\circle*{1}}
  \put(31,24){\circle*{1}}
  \put(20,33){\circle*{1}}
  \put(21,40){\circle*{1}}
  \put(40,33){\circle*{1}}
  \put(41,39){\circle*{1}}
  \put(30,30){\circle*{1}}
  \put(33,44){\circle*{1}}
  \put(-6,2){\makebox(20,10){$(0,0,0)$}}
  \put(-7,18){\makebox(20,10){$(0,0,1)$}}
  \put(42,25){\makebox(10,10){$(2,1,2)$}}
  \put(45,35){\makebox(10,10){$(1,2,2)$}}
  \put(28,42){\makebox(10,10){$(0,2,2)$}}
  \put(7,35){\makebox(10,10){$(0,1,2)$}}
  \put(7,28){\makebox(10,10){$(1,0,2)$}}
  \put(24,22){\makebox(10,10){$(2,0,2)$}}
  \put(28,14){\makebox(20,10)[t]{$(1,1,1)$}}
   \put(32,0){\makebox(10,10)[br]{$x$}}
  \put(42,15){\makebox(10,10){$y$}}
  \put(3,40){\makebox(10,10){$z$}}
  \end{picture}&
  
  \begin{picture}(50,40)(-5,0)
  \put(10,10){\vector(3,-1){28}}
  \put(10,10){\vector(3,1){35}}
  \put(10,10){\vector(0,1){35}}
  \put(10,10){\line(1,1){10}}
  \put(10,10){\line(2,3){11}}
  \put(10,10){\line(3,2){20}}
  \put(10,23){\line(3,-1){10}}
  \put(10,23){\line(3,1){11}}
  \put(21,27){\line(3,-1){10}}
  \put(20,20){\line(3,1){11}}
  \put(10,22.5){\line(1,2){10}}
  \put(21,27){\line(1,1){24}}
  \put(20,20){\line(0,1){23}}
  \put(20,42.5){\line(3,1){25}}
  \put(31.5,24){\line(1,2){13.5}}
  \put(10,10){\circle*{1}}
  \put(10,23){\circle*{1}}
  \put(21,27){\circle*{1}}
  \put(20,20){\circle*{1}}
  \put(31,24){\circle*{1}}
  \put(20,42.5){\circle*{1}}
  \put(45,51){\circle*{1}}
  \put(-6,2){\makebox(20,10){$(0,0,0)$}}
  \put(-7,18){\makebox(20,10){$(0,0,1)$}}
  \put(28,16){\makebox(20,10)[t]{$(1,1,1)$}}
  \put(12,25){\makebox(20,10){$(0,1,1)$}}
  \put(16,12){\makebox(20,10){$(1,0,1)$}}
  \put(32,0){\makebox(10,10)[br]{$x$}}
  \put(42,15){\makebox(10,10){$y$}}
  \put(3,40){\makebox(10,10){$z$}}
  \put(16,42){\makebox(10,10){$(1,0,3)$}}
  \put(47,42){\makebox(10,10){$(1,2,3)$}}
\end{picture}\\
\mbox{(a): $(0,0,1)$ and $(1,1,1)$ of type (1,1)}  &
\mbox{\begin{tabular}{l} (b): $(0,0,1)$ of type (2,0) \\ and $(1,0,1)$ of type (1,0)\end{tabular}}
\end{tabular}
\end{center}
\caption{nonsingular vertex $(0,0,1)$ of types (1,1) and (2,0)}
\label{fig10}
\end{figure}
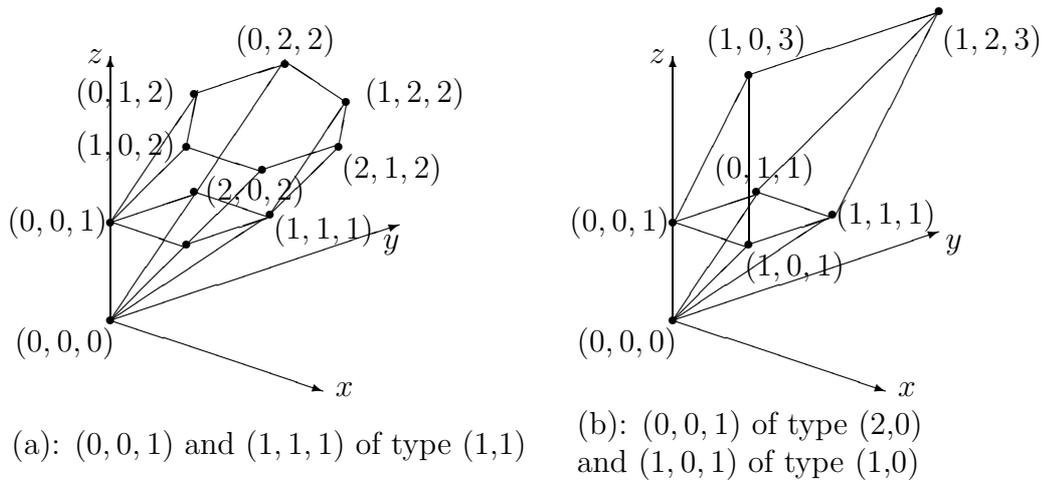

\begin{prop}\label{sect5:p3}
Assume that $\mbox{\rm Int}(Q)\cap M =\emptyset$.
If $Q$ has singular vertices of type $Q_2$, then
$Q$ is normally generated.
\end{prop}

{\it Proof.}
In the Figure~\ref{fig6}, if all four points $(0,0,1), (0,1,1), (1,1,1)$ and $(1,0,1)$ are not vertices,
then all four points $(0,0,2),(0,2,2),(2,2,2)$ and $(2,0,2)$ are vertices contained in one facet,
otherwise the point $(1,1,2)$ is contained in the interior of $Q$. 
In this case, $Q$ is normally generated since $Q$ is isomorphic to the twice of $Q_2$.

If some of four points $(0,0,1), (0,1,1), (1,1,1)$ and $(1,0,1)$ are nonsingular vertices of type (1,1)
and rest are not vertices as in the Figure~\ref{fig10} (a), 
then $Q$ has a facet contained in the plane $\{z=2\}$.
 In this case we cut $Q$ first at the plane $\{z=1\}$ next at the planes $\{y=\pm1\}, \{x+y=3\}$
 and $\{x+y=1\}$ into nonsingular pieces.  Hence, it is normally generated.
 
 If the point $(0,0,1)$ is a singular vertex, then we see $Q$ is normally generated from 
 Proposition~\ref{sect5:p11}.
 Hence, we may assume that an edge connecting two singular vertices of $Q$ has lattice points 
 more than its end points.
 
 We distinguish two cases when $(0,0,1)$ is a nonsingular vertex of type (1,0) or (2,0).
 
 Case (1,0):  The vertex $(0,0,1)$ has two edges connecting with $(1,0,2)$ and $(1,0,1)$, hence,
 the point $(0,1,1)$ is also a vertex of type (1,0) and the other edge from
 $(0,1,1)$ goes to  $(1,2,2)$.   In this case, 
 $Q$ is surrounded by two parallel planes $\{z=x\}$ and $\{z=x+1\}$
 of distance one.
 
 \noindent
 Set $F_{-}:=Q\cap\{z=x\}$ and $F_{+}:=Q\cap \{z=x+1\}$ the parallel facets of $Q$.
 We note that $F_{\pm}$ are nonsingular polygons since the vertices of $Q$ are on 
 one of facets $F_{+}$ and $F_{-}$ and since the line segment of primitive points on the two
 edges from the singular vertex of type $Q_1$ has to be distance two.
 
 \noindent
 Cut off all tetragonal cones with   singular vertices on $F_{-}$ as their apexes
 to obtain an integral polytope $Q'$ with parallel facets
 $F_{+}$ and new $F'\subset F_{-}$ since the distances between singular vertices are
 more than one.
 If $F_{+}$ defines a nonsingular toric surface $Y$, then $F'$ defines a nef divisor on $Y$.
 From Remark~\ref{sect2:rm1}, we see $Q'$ is normally generated.
 
 Case (2,0):  The vertex $(0,0,1)$ has two edges connecting $(1,0,3)$ and $(0,1,1)$,
 hence the point $(0,1,1)$ is also a vertex of type (2,0).  In this case, $Q$ is surrounded by 
 two planes $\{z=x\}$ and $\{z=2x+1\}$.
 If the point $(1,0,1)$ is a nonsingular vertex of type (1,1) or (2,0) or is not vertex, then
 the point $(1,1,2)$ is contained in the interior of $Q$.  Thus $(1,0,1)$ is a vertex of type
 (1,0) hence the point $(1,1,1)$ is also a vertex of type (1,0).  See the Figure~\ref{fig10} (b).
 In this case, $Q$ is bounded by the plane $\{x=1\}$ and the point $(1,3,3)$ is also a vertex.
 By cutting $Q$ at the plane $\{z=1\}$, the rest $Q'$ is a convex hull of 
the line segment $E$ and the tetragon 
 $R=\mbox{Conv}\{(1,0,1),(1,1,0),(1,3,3),(1,0,3)\}$ with two
 edges parallel to $E$  of distance one, hence, 
 we see that $Q'$ is normally generated from Remark~\ref{sect2:rm1}.
\hfill $\Box$

Finally it suffices to consider the case that $Q$ has singular vertex of type $Q_1$.

\begin{figure}[h]
\begin{center}
\setlength{\unitlength}{1mm}
\begin{tabular}{lr}
  \begin{picture}(50,55)(-5,0)
  \put(10,10){\vector(3,-1){28}}
  \put(10,10){\vector(3,1){30}}
  \put(10,10){\vector(0,1){35}}
  \put(10,10){\line(2,1){16}}
  \put(10,10){\line(2,3){24}}
  \put(10,23){\line(3,-1){16}}
  \put(10,23){\line(3,1){11}}
  \put(10,23){\line(1,2){6}}
  \put(10,23){\line(1,3){6}}
  \put(21,27){\line(1,-2){4.5}}
  \put(25.5,18){\line(1,1){15}}
  \put(25.5,18){\line(1,2){6}}
  \put(16,35){\line(3,-1){15}}
  \put(16,35){\line(0,1){6}}
  \put(16,40.5){\line(3,1){18}}
  \put(31,30){\line(3,1){10}}
  \put(41,33){\line(-1,2){7}}
  \put(10,10){\circle*{1}}
  \put(10,23){\circle*{1}}
  \put(21,26.5){\circle*{1}}
  \put(25.5,18){\circle*{1}}
  \put(16,35){\circle*{1}}
  \put(16,40.5){\circle*{1}}
  \put(34,46.5){\circle*{1}}
  \put(31,30){\circle*{1}}
  \put(41,33.5){\circle*{1}}
  \put(3,3){\makebox(10,10){$m_0$}}
  \put(-5,18){\makebox(20,10)[l]{$(0,0,1)$}}
  \put(23,13){\makebox(20,10)[b]{$(2,0,1)$}}
  \put(27,22){\makebox(20,10){$(3,0,2)$}}
  \put(37,31){\makebox(10,10){$(3,1,2)$}}
  \put(37,41){\makebox(10,10){$(0,4,2)$}}
  \put(12,40){\makebox(10,10){$(0,1,2)$}}
  \put(2,29){\makebox(10,10){$(1,0,2)$}}
  \put(32,0){\makebox(10,10)[br]{$x$}}
  \put(37,15){\makebox(10,10){$y$}}
  \put(3,40){\makebox(10,10){$z$}}
\end{picture} &

  \begin{picture}(50,55)(-5,0)
  \put(10,10){\vector(3,-1){28}}
  \put(10,10){\vector(3,1){30}}
  \put(10,10){\vector(0,1){35}}
  \put(10,10){\line(2,1){16}}
  \put(10,10){\line(2,3){11}}
  \put(10,23){\line(3,-1){16}}
  \put(10,23){\line(3,1){11}}
  \put(21,27){\line(1,-2){4.5}}
  \put(10,23){\line(1,2){6}}
   \put(25.5,18){\line(1,1){15}}
  \put(25.5,18){\line(1,2){6}}
   \put(16,35){\line(3,-1){15}}
    \put(31,30){\line(3,1){10}}
   \put(16,35){\line(3,1){20}}
   \put(21,26.5){\line(1,1){15}}
   \put(41,33){\line(-1,2){4.5}}
  \put(10,10){\circle*{1}}
  \put(10,23){\circle*{1}}
  \put(21,26.5){\circle*{1}}
  \put(25.5,18){\circle*{1}}
  \put(16,35){\circle*{1}}
  \put(41,33.5){\circle*{1}}
  \put(31,30){\circle*{1}}
  \put(36,42){\circle*{1}}
  \put(3,3){\makebox(10,10){$m_0$}}
  \put(-5,18){\makebox(20,10)[l]{$(0,0,1)$}}
  \put(23,13){\makebox(20,10)[b]{$(2,0,1)$}}
  \put(9,24){\makebox(20,10){$(0,2,1)$}}
    \put(12,35){\makebox(10,10){$(1,0,2)$}}
    \put(28,22){\makebox(20,10){$(3,0,2)$}}
  \put(42,31){\makebox(10,10){$(3,1,2)$}}
  \put(40,39){\makebox(10,10){$(1,3,2)$}}
  \put(32,0){\makebox(10,10)[br]{$x$}}
  \put(37,15){\makebox(10,10){$y$}}
  \put(3,40){\makebox(10,10){$z$}}
\end{picture}\\
\mbox{(a): $(0,0,1)$ and $(2,0,1)$ of type (1,1)} & 
\mbox{\begin{tabular}{l} (b): $(0,0,1)$ of type (1,0) \\ and $(2,0,1)$ of type (1,1)\end{tabular}}
\end{tabular}
\end{center}
\caption{nonsingular vertex $(0,0,1)$ of type (1,1) and (1,0)}
\label{fig11}
\end{figure}
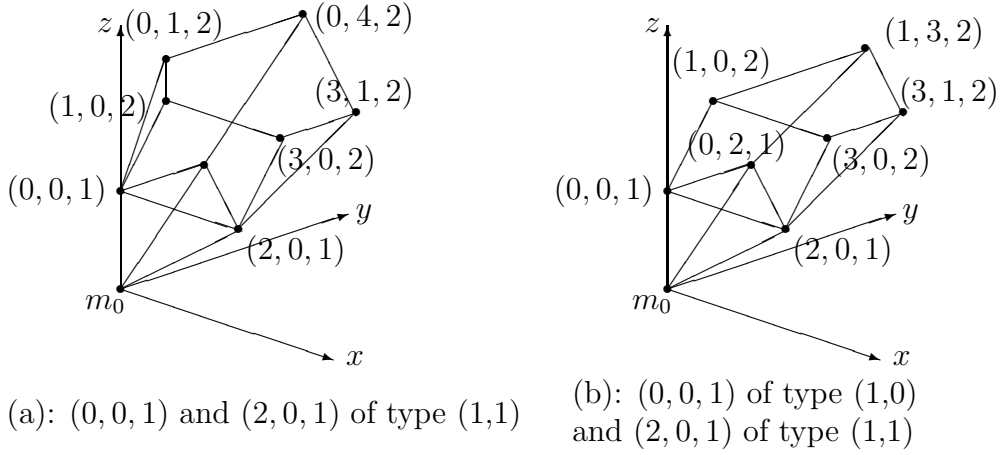

\begin{prop}\label{sect5:p4}
Assume that $\mbox{\rm Int}(Q)\cap M =\emptyset$.
If $Q$ has singular vertices of type $Q_1$, then
$Q$ is normally generated.
\end{prop}

{\it Proof.}
In the Figure~\ref{fig4} (b), if all three points $(0,0,1),(2,0,1)$ and $(0,2,1)$ are not vertices, then
then all points $(0,0,2), (4,0,2)$ and $(0,4,2)$ are contained in one facet, otherwise
the points $(1,1,2),(2,1,2)$ and $(1,2,2)$ are contained in the interior of $Q$.
In this case $Q$ is normally generated since $Q\cong 2Q_1$.
Even if all points $(0,0,1),(2,0,1)$ and $(0,2,1)$ are nonsingular vertices of type (1,1), then
$Q$ is bounded by one facet contained in the plane $\{z=2\}$.  See the Figure~\ref{fig11} (a).
By cutting first at the plane $\{z=1\}$ next at the planes $\{y=x\pm2\}$, the rest is polytope such 
that the bottom is the projective plane and that 
the top facet is nonsingular obtained from the projective plane by blowing ups at two or three points, hence,
it is normally generated.

If the point $(0,0,1)$ is singular of type $Q_2$, then $Q$ is normally generated from 
Proposition~\ref{sect5:p11}.

If all three points $(0,0,1), (0,2,1)$ and $(2,0,1)$ are singular vertices of type $Q_1$, then $Q$ is normally
generated from Proposition~\ref{sect5:p2}.

We assume the point $(0,0,1)$ is a nonsingular vertex of $Q$.

If the point $(0,0,1)$ is a nonsingular vertex of type (2,0), then the other edge from $(0,0,1)$
goes to the point $(1,0,3)$, the point $(0,2,1)$ is also a nonsingular
vertex of type (2,0) and the other edge from $(0,2,1)$ goes to the point $(1,5,3)$. 
In this case, the points $(1,1,2)$ and $(1,2,2)$ are contained in the interior of $Q$.
This case does not occur by assumption.

Next we consider the case when $(0,0,1)$ is a nonsingular vertex of type (1,0). See the Figure~\ref{fig11}
(b).  
In this case, $Q$ has a facet containing in the plane $\{z=x+1\}$, which
 contains the points $(1,1,2)$ and $(1,2,2)$.
If the point $(2,0,1)$ is not vertex, then $Q$ is bounded from above by the plane
$\{z=2\}$ which contains the point $(2,1,2)$.
In this case, after cutting off $Q\cap\{z\le1\}\cong Q_1$ we obtain a polytope with two
parallel facets isomorphic to $2G_0$ and $3G_0$ of distance one, hence, it is normally generated.

If  the point $(2,0,1)$ is a nonsingular vertex,  then it is of type (1,1).  See the Figure~\ref{fig11} (b).
In this case, $Q$ is also bounded by the plane $\{z=2\}$.  The upper facet is a nonsingular tetragon, 
hence, the part $Q\cap\{1\le z\le2\}$ of $Q$ is normally generated as before.

Finally we consider the case that $(0,0,1)$ and $(0,2,1)$ are nonsingular vertices of type (1,0) and $(2,0,1)$
is a singular vertex of type $Q_1$.  In this case, the other two edges from the vertex $(2,0,1)$
go to the points $(2,0,2)$ and $(2,2,2)$.  You may imagine from the Figure~\ref{fig11} (b).

Moreover, if the point $(2,0,2)$ is not vertex, then one edge from $(2,0,1)$ reaches to $(2,0,3)$
and the facet $F_1=Q\cap\{y=0\}$ has the vertex $(2,0,3)$.
If $(2,2,2)$ is not vertex, then the other edge from $(2,0,1)$ reaches to the point $(2,4,3)$.
Set $\tilde Q=\mbox{Conv}\{(0,0,0), (0,0,1), (2,0,1), (2,0,1), (2,0,3), (2,4,3)\}$.
Then we see that $Q$ is contained in $\tilde Q$.

Moreover, if $(2,4,3)$ is a vertex and if $(2,0,3)$ is not vertex of $Q$, then the vertex $(2,4,3)$ has to
be a singular vertex of the facet $Q\cap\{x=2\}$ and it is not $A_1$-singularity.
 Thus we see that if $Q\not=\tilde Q$, then
$Q$ is obtained by cut at the plane $\{z=2\}$  from $\tilde Q$.
If you decompose $\tilde Q$ by cut at the planes $\{z=1\}$ and $\{z=2\}$, then the top part is
isomorphic to $P_{2,3,4}$, hence, it is normally generated.  
The middle part has two facets of distance one, one is isomorphic to $G_{2,3}$ and the other is isomorphic
to $2G_0$.  Since $G_{2,3}$ corresponds to the toric surface obtained by blowing up 
the projective plane at a point, the middle
part is normally generated from Remark~\ref{sect2:rm1}.
Thus we see that $Q$ is normally generated in the case that $(0,0,1)$ is a nonsingular vertex of type (1,0).  
\hfill $\Box$

\section{Proof of Theorem.}\label{sect6}

\begin{prop}\label{sect9:p1}
Let $P$ be a nonsingular convex polytope of dimension three.
If the interior  polytope $Q=\mbox{\rm Conv}\{\mbox{\rm Int}(P) \cap M\}$ of $P$ is of dimension three
without interior lattice points, 
 then $Q$ is normally generated.
\end{prop}

{\it Proof.}
If $Q$ is nonsingular, then it is normally generated from Proposition~\ref{sect2:p2}.
If $Q$ has  singular vertices, then $Q$ is also normally generated
from Proposition~\ref{sect5:p3}.  
\hfill $\Box$

\begin{thm}\label{sect4:tm}
Let $X$ be  a nonsingular projective toric variety
of dimension three and $L$ an ample line bundle on $X$ with
$H^0(X, L\otimes \mathcal{O}_X(2K_X))=0$.  Then $L$ is normally generated.
\end{thm}

{\it Proof.}
Let $P$ be the integral convex polytope of dimension three in $M_{\mathbb{R}}$ corresponding
to the polarized toric variety $(X,L)$.   If $\mbox{Int}(P)\cap M =\emptyset$, then $L$ is
normally generated
 from Proposition~\ref{sect2:p2}.

Consider the case that $\mbox{Int}(P)\cap M \not=\emptyset$.
This implies $\Gamma(X, L\otimes \mathcal{O}_X(K_X))\not=0$.
From Proposition~\ref{sect3:p1}, we have a surjective equivariant morphism $\pi: X\to Y$
 to a nonsingular toric variety $Y$ of dimension three and an ample line bundle $A$ on $Y$
 such that $A\otimes\mathcal{O}_Y(K_Y)$ is generated by global sections and that
 an inclusion $\pi^*(A\otimes\mathcal{O}_Y(K_Y)) \hookrightarrow L\otimes \mathcal{O}_X(K_X)$
 induces the isomorphism of the spaces of global sections.

Set $Q=\mbox{Conv}(\mbox{Int}(P)\cap M)$
the interior polytope of $P$.   The polytope $Q$ corresponds to the globally generated line bundle
$A\otimes\mathcal{O}_Y(K_Y)$ on $Y$.
 If $\dim\ Q\le2$, then  
we can decompose $P$ into a union of $P(F_i)$ with all facets $F_i$ of $P$, and we see that
$L$ is normally generated from
Lemmas~\ref{sect2:l2} and \ref{sect2:l3}.

When $\Gamma(X, L\otimes \mathcal{O}_X(K_X))\not=0$, the vanishing 
$\Gamma(X, L\otimes \mathcal{O}_X(2K_X))=0$ implies the vanishing
$\Gamma(Y, A\otimes\mathcal{O}_Y(2K_Y))=0$
because the inclusion $\pi^*(A\otimes\mathcal{O}_Y(K_Y))\otimes\mathcal{O}_X(K_X)
  \hookrightarrow L\otimes \mathcal{O}_X(2K_X)$
 induces an inclusion $\Gamma(Y, A\otimes\mathcal{O}_Y(2K_Y)) \to
 \Gamma(X, L\otimes\mathcal{O}_X(2K_X))$.
 When $\dim\ Q=3$, thus, the assumption $\Gamma(X, L\otimes\mathcal{O}_X(2K_X))=0$ implies 
 that $\mbox{\rm Int}(Q)\cap M=\emptyset$. 
 This 
  $Q$ is also normally
 generated by Proposition~\ref{sect9:p1}. 
Applying Lemma~\ref{sect2:l2}, we obtain a proof of Theorem~\ref{sect4:tm}.
\hfill $\Box$

\medskip
\section{Application.}\label{sect7}
A nonsingular projective variety $Y$ is called {\it Fano} if its anti-canonical bundle
$\mathcal{O}_Y(-K_Y)$ is ample.

\begin{prop}\label{sect7:p1}
Let $X$ be a nonsingular toric Fano variety of dimension four.  Then
the anti-canonical line bundle of $X$ is normally generated.
\end{prop}

{\it Proof.}
Set $L=\mathcal{O}_X(-K_X)$.  Let $D=\sum_i D_i$ be the divisor consisting all
$T_N$-invariant irreducible divisors on $X$.  Then we have an exact sequence:
\begin{equation}\label{sect5:eq1}
0 \to \mathcal{O}_X \to L \to L_{|D} \to 0.
\end{equation}
Set $P$ the integral convex polytope of dimension four corresponding to the
polarized toric variety $(X, L)$.  
We note that the interior of $P$ contains only one lattice point because 
$L(K_X)\cong \mathcal{O}_X$.
Then the vector space $\Gamma(D, L_{|D}^{\otimes l})$
has a basis $\{{\bold e}(m): m\in \partial(lP)\cap M\}$.   For each $m\in \partial(lP)\cap M$,
we can find $D_i$ so that ${\bold e}(m)
\in \Gamma(L_{|D_i}^{\otimes l})$.
If we could prove the normal generation of $L_{|D_i}$, then we would prove the theorem.

From  Theorem~\ref{sect4:tm}, if each divisor $D_i$ satisfies
$H^0(D_i, L_{|D_i}\otimes \mathcal{O}_{D_i}(2K_{D_i})=0$, then we see that $L_{|D_i}$ and $L$
are normally generated.
We note that $L(-D_i)=\mathcal{O}_X(-K_X-D_i)$ is generated by global sections from 
\cite[Lemma 4.4]{Ms}.
Hence each $D_i$ has globally generated anti-canonical bundle.
By taking a suitable coordinates $(x,y,z,w)$ in $M_{\mathbb{R}}$, we may assume that
$P$ is contained in the half space $\{w\ge0\}$ and that a face $F_i$ of dimension three of $P$
corresponding to $L_{|D_i}$ is $P\cap \{w=0\}$.  Then the globally generated bundle
$L(-D_i)$ corresponds to $P\cap \{w\ge1\}$ and its restriction to $D_i$ does to
the face $G_i:=P\cap \{w=1\}$.
If $H^0(D_i, L_{|D_i}\otimes \mathcal{O}_{D_i}(2K_{D_i})\not=0$, then 
we have an injective homomorphism $\mathcal{O}_{D_i} \to 
L_{|D_i}\otimes \mathcal{O}_{D_i}(2K_{D_i})$.
By tensoring with $\mathcal{O}_{D_i}(-2K_{D_i})$, we have an injection
$\Gamma(D_i,\mathcal{O}_{D_i}(-2K_{D_i}) \to \Gamma(D_i, L_{|D_i}))$, which implies 
$2G_i \subset F_i$ by identifying $M\cap\{w=0\}$ with $M\cap\{w=1\}$.

Let $P':=\mbox{\textup{Conv}}\{(0,0,0,2), 2G_i\}$.  Since $P'\cap\{w=0\}=2G_i \subset
P\cap\{w=0\}=F_i$ and $P'\cap\{w=1\}=G_i = P\cap \{w=1\}=G_i$,
 we have $P'\cap\{w\le1\}\subset P\cap\{w\le1\}$.
 We note $P\cap\{w\le1\}\not=P$ since $(\mbox{Int}\  P)\cap M\not= \emptyset$. 
 
 If $2G_i\not= F_i$, then  $P$ does not contain $(0,0,0,2)$, hence $P$ is
contained in $\{0\le w\le1\}$.  This contradicts to the assumption.
If $2G_i=F_i$, then $P'=P$, hence, $G_i$ is the standard 3-simplex, that is, isomorphic to
$\mbox{\textup{Conv}}\{0, (1,0,0), (0,1,0), (0,0,1)\}$ since $P$ is nonsingular.
In this case, we see that $(X, L) \cong (\mathbb{P}^4, \mathcal{O}(2))$, hence,
$L$ is not the anti-canonical bundle.

From this we see that each divisor $D_i$ satisfies $H^0(D_i, L_{|D_i}\otimes \mathcal{O}_{D_i}(2K_{D_i})=0$.
Hence $L_{|D_i}$ are normally generated for all $i$.  This completes the proof.
\hfill $\Box$

\bigskip
%%%%%%%%%%%% '˜ŽÒŠ'® %%%%%%%%%%%%%
% 'æˆê'˜ŽÒ

\end{document}